\documentclass[12pt]{article}  
\usepackage{amssymb,amsmath,bm,graphicx}
\usepackage{mathrsfs}
\usepackage{constants}
\usepackage{hyperref}
\usepackage{enumerate}
\newcommand{\bsxi}{\boldsymbol{\xi}}

\usepackage{color}
\newcommand{\ncolor}[1]{\textcolor{blue}{#1}}  % Nathakhun's colored text

\begin{document}
\newcommand{\bea}{\begin{eqnarray}}
\newcommand{\ena}{\end{eqnarray}}
\newcommand{\beas}{\begin{eqnarray*}}
\newcommand{\enas}{\end{eqnarray*}}
\newcommand{\beq}{\begin{equation}}
\newcommand{\enq}{\end{equation}}
\def\qed{\hfill \mbox{\rule{0.5em}{0.5em}}}
\newcommand{\bbox}{\hfill $\Box$}
\newcommand{\ignore}[1]{}
\newcommand{\ignorex}[1]{#1}
\newcommand{\wtilde}[1]{\widetilde{#1}}
\newcommand{\qmq}[1]{\quad\mbox{#1}\quad}
\newcommand{\qm}[1]{\quad\mbox{#1}}
\newcommand{\nn}{\nonumber}
\newcommand{\Bvert}{\left\vert\vphantom{\frac{1}{1}}\right.}
\newcommand{\To}{\rightarrow}
\newcommand{\E}{\mathbb{E}}
\newcommand{\Var}{\mathrm{Var}}
\newcommand{\Cov}{\mathrm{Cov}}
\makeatletter
\newsavebox\myboxA
\newsavebox\myboxB
\newlength\mylenA
\newcommand*\xoverline[2][0.70]{%
    \sbox{\myboxA}{$\m@th#2$}%
    \setbox\myboxB\null% Phantom box
    \ht\myboxB=\ht\myboxA%
    \dp\myboxB=\dp\myboxA%
    \wd\myboxB=#1\wd\myboxA% Scale phantom
    \sbox\myboxB{$\m@th\overline{\copy\myboxB}$}%  Overlined phantom
    \setlength\mylenA{\the\wd\myboxA}%   calc width diff
    \addtolength\mylenA{-\the\wd\myboxB}%
    \ifdim\wd\myboxB<\wd\myboxA%
       \rlap{\hskip 0.5\mylenA\usebox\myboxB}{\usebox\myboxA}%
    \else
        \hskip -0.5\mylenA\rlap{\usebox\myboxA}{\hskip 0.5\mylenA\usebox\myboxB}%
    \fi}
\makeatother

\newtheorem{theorem}{Theorem}[section]
\newtheorem{corollary}[theorem]{Corollary}
\newtheorem{conjecture}[theorem]{Conjecture}
\newtheorem{proposition}[theorem]{Proposition}
\newtheorem{lemma}[theorem]{Lemma}
\newtheorem{definition}[theorem]{Definition}
\newtheorem{example}[theorem]{Example}
\newtheorem{remark}[theorem]{Remark}
\newtheorem{case}{Case}[section]
\newtheorem{condition}{Condition}[section]
\newcommand{\proof}{\noindent {\it Proof:} }

\title{{\bf\Large Stein's method for positively associated random variables with applications to the Ising and voter models, bond percolation, and contact process}}
\author{Larry Goldstein\thanks{Research partially supported by NSA-H98230-15-1-0250.}  \ and Nathakhun Wiroonsri\\University of Southern California}
\maketitle
\date{}

\begin{abstract}   
We provide non-asymptotic $L^1$ bounds to the normal for four well-known models in statistical physics and particle systems in $\mathbb{Z}^d$; the ferromagnetic nearest-neighbor Ising model, the supercritical bond percolation model, the voter model and the contact process. In the Ising model, we obtain an $L^1$ distance bound between the total magnetization and the normal distribution at any temperature when the magnetic moment parameter is nonzero, and when the inverse temperature is below critical and the magnetic moment parameter is zero. In the percolation model we obtain such a bound for the total number of points in a finite region belonging to an infinite cluster in dimensions $d \ge 2$, in the voter model for the occupation time of the origin in dimensions $d \ge 7$,  and for finite time integrals of non-constant increasing cylindrical functions evaluated on the one dimensional supercritical contact process started in its unique invariant distribution.

The tool developed for these purposes is a version of Stein's method adapted to positively associated random variables.
%$\{ X_{{\bf j}}:{\bf j} \in \mathbb{R}^d \}$. 
In one dimension, letting $\bsxi=(\xi_1,\ldots,\xi_m)$ be a positively associated mean zero random vector with components that obey the bound $|\xi_i| \le B, i=1,\ldots,m$, and whose sum $W = \sum_{i=1}^m \xi_i$ has variance 1, it holds that
\beas
d_1 \big({\cal L}(W),{\cal L}(Z) \big) \le 5B  + \sqrt{\frac{8}{\pi}}\sum_{i \not = j} \E[\xi_i \xi_j]
\enas
where $Z$ has the standard normal distribution and $d_1(\cdot,\cdot)$ is the $L^1$ metric. Our methods apply in the multidimensional case with the $L^1$ metric replaced by a smooth function metric.
\end{abstract}

\section{Introduction}
We provide non-asymptotic $L^1$ bounds to the normal for four well-known models in statistical physics and particle systems in $\mathbb{Z}^d$; the ferromagnetic nearest-neighbor Ising model, supercritical bond percolation, the voter model and the contact process. We recall that the $L^1$, or Wasserstein, distance between the distributions ${\cal L}(X)$ and ${\cal L}(Y)$ of real valued random variables $X$ and $Y$ is given by
\bea \label{def.d1.integral}
d_1 \big({\cal L}(X),{\cal L}(Y)\big) = \int_{-\infty}^\infty |P(X \le t) - P(Y \le t)| dt.
\ena

Taking the Ising model first, questions regarding the distribution of total magnetization in statistical mechanical models have attracted attention for quite some time. To fix terms, with $\mathbb{Z}$ the set of integers and $k \in \mathbb{Z}$ let $\mathbb{N}_k=[k,\infty) \cap \mathbb{Z}$, and for $d$ a positive integer let $\mathbb{Z}^d$ denote the $d$-dimensional lattice composed of vectors ${\bf i}=(i_1,\ldots,i_d)$ having integer components, and let $\Lambda$ be a symmetric hypercube, that is,
\beas 
\Lambda = \{ {\bf i} \in \mathbb{Z}^d: |{\bf i}|_{\infty} \le k \} \text{ \ for some \ } k \in \mathbb{N}_0,
\enas
where 
$|{\bf i}|_{\infty} = \sup_{1 \le j \le d}|i_j|$, the supremum norm. We say $\omega =\{ \omega_{\bf i}: {\bf i} \in \Lambda\}$ is a configuration on $\Lambda$ when $\omega_{{\bf i}}$, denoting the magnetic spin at site ${\bf i}$, takes values in the set $\{-1,+1\}$ for each ${\bf i} \in \Lambda$. 

We define a probability measure on configurations on $\Lambda$, corresponding to the nearest-neighbor ferromagnetic Ising model. First, fix a configuration $\widetilde{\omega}=\{ \widetilde{\omega}_{\bf i}: {\bf i} \in \mathbb{Z}^d\}$ on all of $\mathbb{Z}^d$, and the interaction strength $J \ge 0$. Now define the energy function, or Hamiltonian $H_{\Lambda,h,\widetilde{\omega}}$, by
\bea \label{hamil}
H_{\Lambda,h,\widetilde{\omega}}(\omega) = - J\left(\sum_{\substack{ {\bf i},{\bf j} \in \Lambda \\ |{\bf i}-{\bf j}|_{1}=1}}\omega_{\bf i} \omega_{\bf j}+\sum_{\substack{ {\bf i} \in \Lambda,{\bf j} \notin \Lambda \\ |{\bf i}-{\bf j}|_{1}=1}}\omega_{\bf i} \widetilde{\omega}_{\bf j}\right) - h \sum_{{\bf i} \in \Lambda} \omega_{{\bf i}},
\ena
where $h$ is the magnetic moment parameter and $|{\bf x}|_1=\sum_{i=1}^d |x_i|$ is the $L^1$ norm of ${\bf x} \in \mathbb{R}^d$. 
Given an inverse temperature $\beta>0$, the finite volume Ising probability measure on $\Lambda$ is given by
\beas
\mu_{\Lambda,\beta,h,\widetilde{\omega}}(\omega) = \exp (-\beta H_{\Lambda,h,\widetilde{\omega}}(\omega))/Z_{\Lambda, \beta, h, \widetilde{\omega}},
\enas
for a normalizing constant $Z_{\Lambda, \beta, h, \widetilde{\omega}}$. The cases where $\widetilde{\omega}_{{\bf j}}$ takes the value $0,+1$ and $-1$ for all ${\bf j} \in \mathbb{Z}^d$ are typically referred to as {\it free, positive, and negative boundary conditions}, respectively, see \cite{ELL06}, for example.

Now let $\{ \Lambda_n:n\in \mathbb{N}_1\}$ be any increasing sequence of symmetric hypercubes of $\mathbb{Z}^d$ whose union is $\mathbb{Z}^d$. For fixed $\beta$ and $h$, we say that $\mu_{\Lambda_n,\beta,h,\widetilde{\omega}}$ converges \textit{weakly} to $\mu_{\beta,h,\widetilde{\omega}}$ as $n \rightarrow \infty$ if 
\beas
\lim_{n\rightarrow \infty} \int f d\mu_{\Lambda_n,\beta,h,\widetilde{\omega}} = \int f d\mu_{\beta,h,\widetilde{\omega}}
\enas
for all local functions $f:\mathbb{R}^d \rightarrow \mathbb{R}$, where we recall that a function $f$ is said to be \textit{local} if it depends only on $\Lambda_n$ for some $n$. Each such weak limit is a probability measure on $\{-1,1 \}^{\mathbb{Z}^d}$. The set of infinite volume Gibbs measures, denoted ${\cal G}_{\beta,h}$, is the closed convex hull of all weak limits of $\mu_{\Lambda_n,\beta,h,\widetilde{\omega}}$ for $\widetilde{\omega}$ any configuration on $\mathbb{Z}^d$.

The structure of $\mathcal{G}_{\beta,h}$ depends on parameters $\beta$ and $h$. In particular, by \cite{ELL06}, there exists a critical inverse temperature $\beta_c$ such that the limit measures depend, or do not depend, on the boundary conditions when $\beta > \beta_c$ or $\beta < \beta_c$, respectively. More precisely, for $\beta >0, h\ne 0$ and $0< \beta < \beta_c,h=0$, the finite volume measures $\mu_{\Lambda_n,\beta,h,\widetilde{\omega}}$ have a unique weak limit for any choice of boundary conditions $\widetilde{\omega}$. Thus, in these cases, $\mathcal{G}_{\beta,h}$ consists of a unique measure, which we denote $\mu_{\beta,h}$. In the case $\beta > \beta_c, h = 0$, the set $\mathcal{G}_{\beta,0}$ contains two \textit{pure phases}, $\mu^{+}_{\beta,0}$ and $\mu^{-}_{\beta,0}$, which are the infinite volume measures arising as weak limits under positive and negative boundary conditions respectively. Hence, by definition, $\mathcal{G}_{\beta,0}$ also contains the convex combinations $\mu^{(\alpha)}_{\beta,0} = \alpha \mu^{+}_{\beta,0} + (1-\alpha)\mu^{-}_{\beta,0}$, $0 < \alpha < 1$. These measures are called \textit{mixed phases}. For $d=2$, $\mathcal{G}_{\beta,0}$ consists of only pure and mixed phases. For $d \ge 3$, the set $\mathcal{G}_{\beta,0}$ contains nontranslation invariant measures for any sufficiently large $\beta$, which are not considered in this paper. For additional detail about finite and infinite volume Gibbs measures, see \cite{ELL06}.

In this work, we focus on the Ising models on $\mathbb{Z}^d$ in the cases where the weak limit is unique, precisely, for measures in the set
\bea \label{Ismeasures1}
\mathcal{M}_1 = \{\mu_{\beta,h} :(\beta,h) \in \mathbb{R}^+ \times \mathbb{R} \setminus \{0 \} \cup (0,\beta_c)\times \{0 \} \}.
\ena
For these measures, we obtain $L^1$ bounds to the normal for the distribution of the total magnetization 
\bea \label{mkdef}
M_{\bf k} = \sum_{{\bf i} \in B_{{\bf k}}^n}\omega_{{\bf i}}
\ena
over the finite block $B_{{\bf k}}^n \subset \mathbb{Z}^d$ with side length $n$ `anchored' at ${\bf k} \in \mathbb{Z}^d$, defined in \eqref{blockdef}.
See Remark \ref{re:AD15.etc} for a discussion regarding the applicability of our results to the
pure or mixed phases when the inverse temperature is above critical, that is, for measures in the set
\bea \label{Ismeasures2}
\mathcal{M}_2 = \begin{Bmatrix} \mu^{+}_{\beta,0},\mu^{-}_{\beta,0} &:&\beta > \beta_c,h=0 \\ \mu^{(\alpha)}_{\beta,0} &:&0 < \alpha < 1,\beta > \beta_c,h=0 \end{Bmatrix}.
\ena

The second model we consider is bond percolation on $\mathbb{Z}^d$ for $d \ge 2$ (see \cite{Gri99}, for example). In graphical terms, nearest neighbor bonds are the edges between vertices ${\bf x}$ and ${\bf y}$ in $\mathbb{Z}^d$ satisfying $|{\bf x}-{\bf y}|_{1}=1$; we denote the collection of all such bonds by
$\mathbb{E}^d$. For a given $\theta \in [0,1]$, we declare each bond in $\mathbb{E}^d$ to be {\it open} with probability $\theta$ and {\it closed} otherwise, independently of all other bonds. More formally, we take ${\bf G} = \{ 0,1 \}^{\mathbb{E}^d}$, elements of which are represented as $\boldsymbol{g} = \left\{ g({\bf e}):{\bf e} \in \mathbb{E}^d\right\}$, and called configurations. The value $g({\bf e})=0$ corresponds to the bond ${\bf e}$ being closed and $g({\bf e})=1$ corresponds to ${\bf e}$ being open. 

Given a bond configuration $\boldsymbol{g}$ and ${\bf x} \in \mathbb{Z}^d$, 
let $\mathcal{C}({\bf x})$ be the set of vertices connected to ${\bf x}$ through open bonds, and let $|\mathcal{C}({\bf x})|$ denote the number of vertices in $\mathcal{C}({\bf x})$, which may be infinite. Since the probability $\theta$ that a bond is open is the same for every ${\bf e} \in \mathbb{E}^d$ and each bond is independent of all others, $\mathcal{C}({\bf x})$ has the same distribution for all ${\bf x} \in \mathbb{R}^d$. We define
\beas %\label{infclusprob}
\rho(\theta) = P(|\mathcal{C}({\bf 0})|= \infty) \qmq{and}  \theta_c = \sup\{ \theta: \rho(\theta)=0 \},
\enas
respectively representing the probability that for a given connectivity $\theta$ a vertex belongs to an infinite cluster, and its threshold, the \textit{critical probability}. 
Our interest is in obtaining $L^1$ bounds to the normal, in the supercritical case $\theta>\theta_c$, for the distribution of the total number of points in a finite block that belong to an infinite cluster, specifically for
\bea \label{ukdef}
U_{\bf k} = \sum_{{\bf x} \in B_{{\bf k}}^n}\mathbf{1}_{ \{ |\mathcal{C}({\bf x})|=\infty \} }
\ena
where $B_{{\bf k}}^n \subset \mathbb{Z}^d$ is defined in \eqref{blockdef}.

The third application we consider is the $d$-dimensional voter model introduced in \cite{CS73}, a Markov process $\{ \eta_t:t \ge 0 \}$ taking values in $\{ 0,1 \}^{\mathbb{Z}^d}$, and one of the simplest interacting particle systems, see also \cite{cox83}. With a precise formulation deferred to Section \ref{app3}, the process is specified by
taking the initial distribution $\{ \eta_0(\mathbf{x}): \mathbf{x} \in \mathbb{Z}^d \}$ to be a family of independent, identically distributed Bernoulli random variables with parameter $\theta \in (0,1)$, and having transitions
\beas
\eta_t(\mathbf{x}) \rightarrow 1-\eta_t(\mathbf{x}) \text{ \ at rate \ } \frac{1}{2d}\big|\{ \mathbf{y}:|\mathbf{y}-\mathbf{x}|_{1}=1 \text{ \ and \ } \eta_t(\mathbf{x}) \neq \eta_t(\mathbf{y})\}\big|.
\enas
That is, a site ${\bf x} \in \mathbb{Z}^d$ changes its state with rate equal to the fraction of its $2d$ nearest neighbors that are in the opposite state. We study the occupation time
\bea \label{occudef}
T_s^t = \int_{s}^{s+t} \eta_u(\mathbf{0})du \qm{for $t > 0$, $s \ge 0$,}
\ena
that is, the amount of time in $(s,s+t]$ that the origin ${\bf 0}$ spends in state 1.

The last model we consider is the one dimensional contact process $\{\zeta^{\nu_\lambda}(t):t\ge 0\}$ introduced in \cite{Har74}, a continuous time Markov process with state space $\mathscr{P}(\mathbb{Z})$, the set of all subsets of $\mathbb{Z}$. The subset $\zeta^{\nu_\lambda}(t)$ models the collection of `infected' individuals at time $t$, where each infected site infects healthy neighbors at rate $\lambda$, and  infected sites recover at rate $1$. The distribution of $\zeta^{\nu_\lambda}(0)$ is $\nu_{\lambda}$,  the unique invariant measure defined in \cite{DG83} as the limiting distribution of $\zeta^{\mathbb{Z}}(t)$,
\beas
\zeta^{\mathbb{Z}}(t)  \overset{d}{\to} \nu_\lambda \text{ \ as \ } t \rightarrow \infty,
\enas
where $\zeta^{\mathbb{Z}}(t)$ has initial distribution putting mass one on all of $\mathbb{Z}$.

Define
\bea \label{contactcritical}
\lambda_* = \sup \{\lambda >0: \zeta^{\mathbb{Z}}(t) \rightarrow \delta_{\emptyset} \text{ weakly as } t \rightarrow \infty \}.
\ena 
The contact process is said to be supercritical when $\lambda>\lambda_*$, and, by \cite{DG83}, in this case $\nu_\lambda$ is the nontrivial, unique invariant measure of the process; $\nu_{\lambda} = \delta_{\emptyset}$ when $\lambda<\lambda_*$.

Recall that a function $f:\mathscr{P}(\mathbb{Z}) \rightarrow \mathbb{R}$ is said to be \textit{cylindrical} if it depends on only finitely many sites and \textit{increasing} if $f(A) \le f(B)$ whenever $A \subset B \subset \mathbb{Z}$. In the supercritical case, we study the cumulative value of a cylindrical $f$ evaluated on the process over the interval $(s,s+t]$,
\bea \label{contactdef}
D_{s,f}^t = \int_s^{s+t} f(\zeta^{\nu_\lambda}(u)) du \text{ \ for \ } s\ge 0,t > 0.
\ena
For instance, for $B$ a finite subset of $\mathbb{Z}$, letting $I_{\xoverline{B}}(\cdot)$ be the cylindrical, increasing indicator function given by 
\bea \label{indicatordef}
I_{\xoverline{B}}(\eta)={\bf 1}(\eta \cap B \ne \emptyset) \qmq{for $\eta \subset \mathbb{Z}$,}
\ena
the value $D_{s,I_{\bar{B}}}^t$ yields the amount of time in $(s,s+t]$ that $B$ contains at least one infected site.

The quantities \eqref{mkdef}, \eqref{ukdef}, \eqref{occudef} and \eqref{contactdef} that we consider in these four models are obtained by summing or integrating positively associated variables, and we obtain our results by developing an $L^1$ version of Stein's method for such instances. 
We recall that a random vector $\bsxi=(\xi_1,\ldots,\xi_m) \in \mathbb{R}^m$ is said to be \textit{positively associated} whenever
\bea \label{fkg:def}
\Cov\big(\psi(\bsxi),\phi(\bsxi)\big) \ge 0
\ena
for all real valued coordinate-wise nondecreasing functions $\psi$ and $\phi$ on $\mathbb{R}^m$ such that 
$\psi(\bsxi)$ and $\phi(\bsxi)$ posses second moments. In general, a collection $\{\xi_\alpha: \alpha \in I\}$ of real valued random variables indexed by a set $I$ is said to be positively associated if all finite subcollections are positively associated. Positive association was introduced in \cite{EPW67} and has been found frequently in probabilistic models in several areas, especially statistical physics. In some literature \eqref{fkg:def} is termed the FKG-inequality, see \cite{new80} for example.

For the Ising model with measure in $\mathcal{M}_1 \cup \mathcal{M}_2$ as given in \eqref{Ismeasures1} and \eqref{Ismeasures2}, it is well-known that the variables making up the configuration $\{ \omega_{{\bf i}}:{\bf i} \in \mathbb{Z}^d \}$ are positively associated whenever $J \ge 0$. For the percolation model, it is also well-known that $\{ \mathbf{1}_{ \{|\mathcal{C}({\bf x})| = \infty \} }:{\bf x} \in \mathbb{Z}^d \}$ is positively associated (see \cite{new80} for both models). For the voter model, taking $t > 0$, $m \in \mathbb{N}_1$ and
\bea \label{occudef2}
X_{s,i}^t= \int_{s+(i-1)t/m}^{s+it/m} \eta_u(\mathbf{0})du \text{ \ for \ } i \in [m], \qmq{we have} T_s^t = \sum_{i=1}^m X_{s,i}^t,
\ena
where $[m]=\{ 1,2,\ldots,m \}$. 
It was demonstrated in \cite{cox84} that the family $\{ X_{s,i}^t: i \in [m] \}$ is positively associated. For the contact process, similarly, with $t>0$, $m \in \mathbb{N}_1$ and
\bea \label{yidef}
Y_{s,i}^{t,m} = \int_{s+(i-1)t/m}^{s+it/m} f(\zeta^{\nu_\lambda}(u))du \text{ \ for \ } i \in [m], \qmq{we have} D_{s,f}^t = \sum_{i=1}^m Y_{s,i}^{t,m}.
\ena 
We show in Lemma \ref{PAcontact}, using a simple consequence from arguments in the proof of Lemma 1 of \cite{Sch86}, that the family $\left\{ Y_{s,i}^{t,m} : i \in [m]\right\}$ is positively associated when $f$ is increasing. Hence in each of the models we consider the quantity we study can be represented as the sum of positively associated random variables.

Over the last few decades, several papers established central limit theorems and rates of convergence for sums of positively associated random variables under various assumptions. To state some of these results, recall that a random field $\{X_{{\bf j}}: {\bf j} \in \mathbb{Z}^d\}$ is called {\em second order stationary} when $EX_{{\bf j}}^2<\infty$ for all ${\bf j} \in \mathbb{Z}^d$ and the covariance $\Cov\big(X_{{\bf i}},X_{{\bf j}}\big)$ depends only on ${\bf j}-{\bf i}$, that is, $\Cov \big(X_{{\bf i}},X_{{\bf j}} \big) = R({\bf j}-{\bf i})$ for all ${\bf i},{\bf j}\in \mathbb{R}^d$, with $R(\cdot)$ necessarily given by
\bea \label{covdef}
R({\bf k})=\Cov \big(X_{{\bf 0}},X_{{\bf k}} \big).
\ena
We say the field is \textit{strictly stationary} if for all $m \in \mathbb{N}_1$ and ${\bf k},{\bf j}_1,\ldots,{\bf j}_m \in \mathbb{Z}^d$, the vector $\left(X_{{\bf j}_1}, \ldots , X_{{\bf j}_m}\right)$ has the same distribution as $\left(X_{{\bf j}_1+{\bf k}}, \ldots , X_{{\bf j}_m+{\bf k}}\right)$. From the definitions, when second moments exist, strict stationarity implies second order stationarity, and if a second order stationary field is positively associated, then $R({\bf k}) \ge 0$ for all ${\bf k} \in \mathbb{Z}^d$.

We let ${\bf 1} \in \mathbb{Z}^d$ denote the vector with all components $1$, and write inequalities such as ${\bf a} < {\bf b}$ for vectors ${\bf a} , {\bf b} \in \mathbb{R}^d$  when they hold componentwise.
For ${\bf k} \in \mathbb{Z}^d, n \in \mathbb{N}_1$ and a random field $\{X_{{\bf j}}: {\bf j} \in \mathbb{Z}^d\}$, define the `block sum' variables, over a block with side length $n$, by 
\bea \label{blockdef}
S_{{\bf k}}^n =  \sum_{{\bf j} \in B_{{\bf k}}^n} X_{{\bf j}} \qmq{where}
B_{{\bf k}}^n = \left\{ {\bf j} \in \mathbb{Z}^d: {\bf k} \le {\bf j} < {\bf k}+n{\bf 1} \right\}.
\ena
Note that $B_{\bf k}^n=B_{\bf 0}^n+{\bf k}$.

For a second order stationary field with $R(\cdot)$ given by \eqref{covdef}, we have
\bea \label{andef}
\Var \big(S_{{\bf k}}^n \big)  = \sum_{{\bf i},{\bf j} \in B_{{\bf k}}^n} \Cov \big(X_{{\bf i}},X_{{\bf j}} \big) = n^d A_n \qmq{where} A_n = \frac{1}{n^d} \sum_{{\bf i},{\bf j} \in B_{{\bf 1}}^n} R({\bf i}-{\bf j}).
\ena
With
\bea \label{def:A}
A = \sum_{{\bf k} \in \mathbb{Z}^d} R({\bf k}),
\ena 
Lemma 4 of \cite{new80} shows
$R({\bf 0}) < \infty$, $R({\bf k}) \ge 0$ and $A<\infty$ imply $\lim_{n \rightarrow \infty} A_n=A$. As $A_n$ or $A$ equals zero only if the field is trivial, we assume without loss of generality in the following that $A_n>0$ for all $n \in \mathbb{N}_1$ and $A>0$.

Theorem \ref{newmantheorem} of \cite{new80} shows that the collection of properly standardized block sums $S_{{\bf k}}^n$ of a strictly stationary, positively associated random field with covariance given by \eqref{covdef} converges jointly to independent normal variables when $A<\infty$. 

\begin{theorem}[\cite{new80}] \label{newmantheorem}
	Let $\{X_{{\bf j}}: {\bf j} \in \mathbb{Z}^d\}$ be a positively associated strictly stationary random field with finite second moments and covariance $R({\bf k})$ given by \eqref{covdef} that satisfies $A<\infty$ where $A$ is given in \eqref{def:A}. Then the standardized block sums
	\bea \label{eq:SnkNew}
	\left\{\frac{S_{n{\bf k}}^n - \E S_{n{\bf k}}^n}{n^{d/2}}:{\bf k}\in \mathbb{Z}^d \right\}  
	\ena
	converge to independent mean zero normal distributions with variance $A$ as $n \rightarrow \infty$. That is, the expectations of bounded continuous functions of \eqref{eq:SnkNew} that depend on only finitely many coordinates converge to the expectation of that same number of independent mean zero, variance $A$ normal variables. 
\end{theorem}

We note that the blocks $B_{n{\bf k}}^n$ over which sums are taken in Theorem \ref{newmantheorem} are disjoint, and that the limiting distribution has independent coordinates. As an application of our Theorem \ref{stein:fkg2}, in Theorem \ref{multiising2} we obtain bounds for the sums over overlapping blocks to an approximating dependent multivariate normal distribution. Rectangular blocks, that is, blocks with varying side lengths, can also be accommodated by our methods, at the cost of some additional complexity in our computations.

The stationarity assumption was relaxed in \cite{cox84}, where it was shown that the block sums of a positively associated random field $\{X_{{\bf j}}: {\bf j} \in \mathbb{Z}^d\}$ converge to the normal whenever $X_{{\bf j}}$ has finite third moment for every ${\bf j} \in \mathbb{Z}^d$ and $u(n)$ converges to zero as $n \rightarrow \infty$ where 
\bea \label{undef}
u(n) = \sup_{{\bf j} \in \mathbb{Z}^d}\sum_{{\bf i} \in \mathbb{Z}^d : |{\bf j}-{\bf i}|_{\infty} \ge n}\Cov \big(X_{\bf i},X_{\bf j} \big).
\ena

The rate of convergence for sums of positively associated variables was first obtained in \cite{Bir88}, which achieved the rate $\log (n) / \sqrt{n}$ in the Kolmogorov metric in the one dimensional case, assuming uniformly finite $3+\epsilon$-moments for some $\epsilon > 0$, and that $u(n)$ decays at an exponential rate. Theorem \ref{Bultheorem} of \cite{Bul95} generalizes the results in \cite{Bir88}. In the following $C$ will denote a constant whose value may change from line to line. Recall that the Kolmogorov, or $L^\infty$ distance between distributions ${\cal L}(X)$ and ${\cal L}(Y)$ is given by 
\beas %\label{def.dk.integral}
d_K \big({\cal L}(X),{\cal L}(Y)\big) = \sup_{t \in \mathbb{R}} |P(X \le t) - P(Y \le t)|.
\enas

\begin{theorem}[\cite{Bul95}] \label{Bultheorem} 
Let $\epsilon>0$ and $\{X_{{\bf j}}: {\bf j} \in \mathbb{Z}^d\}$ be a mean zero, positively associated random field whose elements have uniformly bounded $3+\epsilon$ moments. Assume that there exists $\lambda > 0$ such that $u(n)$ as defined in \eqref{undef} satisfies
\beas %\label{eq:Bulun}
u(n) \le \kappa_0e^{-\lambda n} \text{ \ \ for some \ \ } \kappa_0 >0 .
\enas
Then for any finite subset $V \subset \mathbb{Z}^d$ and 
\beas %\label{svdef}
S(V) = \sum_{{\bf j} \in V} X_{{\bf j}} \text{ \ with  \ } \sigma^2(V) = {\rm Var}(S(V)),
\enas
one has
\beas
d_K\big({\cal L}(S(V)/\sigma(V)),{\cal L}(Z)\big) \le C|V|\sigma(V)^{-3}(\log(|V|+1))^d 
\enas
where $Z$ is a standard normal random variable, $|V|$ denotes the size of $V$ and $C>0$ is a constant depending only on $\epsilon$, $d$, $\kappa_0$ and $\lambda$.
\end{theorem}

In the present work we first provide bounds in the $L^1$ distance to the normal for a sum of elements of a positively associated random field in $\mathbb{Z}^d$, under the condition that covariances decay at an exponential rate. We apply these bounds to the total magnetization of the $d$-dimensional Ising model as defined in \eqref{mkdef} and the total number of points belonging to an infinite cluster in the bond percolation model as defined in \eqref{ukdef}. 

We also consider the occupation time of the voter model, where the assumption that the covariance decays exponentially is not satisfied. Without using positive association, \cite{cox83} showed that the occupation time $T_0^t$ of the voter model as defined in \eqref{occudef} satisfies the CLT for $d\ge 2$. In later work, \cite{cox84} used the fact that $T_0^t$ is the sum of positively associated random variables and that $u(n)$ as defined in (\ref{undef}) converges to zero as $n \rightarrow \infty$ to prove the CLT for $d \ge 5$. However, error bounds and the rate of convergence appear yet unaddressed in the literature. Bounds to the limit in the $L^1$ metric for $d \ge 7$ are obtained using the technique provided in Theorem \ref{stein:fkg} below.

Next, with definitions as in the Introduction, we provide $L^1$ bounds for the quantity $D_{s,f}^t$ in \eqref{contactdef}, an integral over the time interval $(s,s+t]$ of a non-constant increasing cylindrical function $f$ evaluated on the supercritical one dimensional contact process $\zeta^{\nu_{\lambda}}(t), t \ge 0$ with initial state having distribution $\nu_\lambda$. Our argument exploits the exponential decay of the covariance between $f(\zeta^{\nu_{\lambda}}(r))$ and $f(\zeta^{\nu_{\lambda}}(s))$ for any cylindrical $f$. \cite{Sch86} applied the result in \cite{new83} along with the fact that $D_{0,f}^t$ is the sum of positively associated random variables to prove a CLT for $D_{0,f}^t$. Error bounds in Kolmogorov distance, without explicit constants, can be obtained by applying Theorem \ref{Bultheorem} for our applications to the Ising model, bond percolation and the contact process, contingent upon showing that $u(n)$ in \eqref{undef} decays exponentially.

We obtain $L^1$ results on the four models considered by developing a version of Stein's method for sums of bounded positively associated random variables. Stein's method was first introduced by Charles Stein in his seminal paper \cite{Stein72} and has become one of the most powerful methods to prove convergence in distribution. Its main advantages are that it provides non asymptotic bounds on the distance between distributions, and that it can handle a number of situations involving dependence. For more detail about the method in general, see the text \cite{Che11} and the introductory notes \cite{Ross11}. 

Stein's method has been used previously in several papers in statistical physics. In \cite{ER08}, a version of Stein's method using size bias couplings was developed for a class of discrete Gibbs measures, that is, probability measures on $\mathbb{N}_0$ proportional to $e^{V(k)}$ for some function $V:\mathbb{N}_0 \rightarrow \mathbb{R}$. This work provided bounds in the total variation distance between the distribution of certain sums of strongly correlated random variables and discrete Gibbs measures, with applications to interacting particle systems. Stein's method of exchangeable pairs was applied to the classical Curie-Weiss model for high temperatures and at the critical temperature in \cite{EL10}, where optimal Kolmogorov distance bounds were produced. The same result at the critical temperature was also obtained in \cite{CS11}.  These results were extended to the Curie-Weiss-Potts model in \cite{EM15}. 

Our main result in the one dimensional case is the following.
\begin{theorem}
\label{stein:fkg}
Let $\bsxi=(\xi_1,\ldots,\xi_m)$ be a positively associated mean zero random vector with components obeying the bound $|\xi_i| \le B$, and whose sum $W = \sum_{i=1}^m \xi_i$ has variance 1. Let $Z$ be a standard normal random variable. Then
\bea \label{stein:fkgbound}
d_1\big({\cal L}(W),{\cal L}(Z)\big) \le  5B + \sqrt{\frac{8}{\pi}}\sum_{i \not = j} \sigma_{ij} \qmq{where} \sigma_{ij}=\E[\xi_i \xi_j].
\ena
\end{theorem}

We also consider a version of Theorem \ref{stein:fkg} adapted to the multidimensional case, with the $L^1$ metric replaced by a smooth functions metric, following the development of Chapter 12 of \cite{Che11}. 
For a real valued function $\varphi(u)$ defined on the domain $\mathcal{D}$, let $|\varphi|_{\infty}= \sup_{x \in \mathcal{D}} |\varphi(x)|$. We include in this definition the $|\cdot|_\infty$ norm of vectors and matrices, for instance, by considering them as real valued functions of their indices. 

For $m \in \mathbb{N}_0$, let $L^{\infty}_m(\mathbb{R}^p)$ be the collection of all functions $h:\mathbb{R}^p \rightarrow \mathbb{R}$ such that for all $\mathbf{k} = (k_1,\ldots,k_p) \in \mathbb{N}_0^p$ with $|{\bf k}|_1 \le m$, the partial derivative 
\beas
h^{(\mathbf{k})}(\mathbf{x}) = \frac{\partial^{|\mathbf{k}|_1}h}{\partial^{k_1}x_1\cdots\partial^{k_p}x_p} 
\enas
exists, and 
\beas
|h|_{L^{\infty}_m(\mathbb{R}^p)}: = \max_{0 \le |\mathbf{k}|_1 \le m}|h^{(\mathbf{k})}|_{\infty} \text{ \ \ is finite.}
\enas
Let
\beas
\mathcal{H}_{m,\infty,p} = \{ h\in L^{\infty}_m(\mathbb{R}^p):|h|_{L^{\infty}_m(\mathbb{R}^p)} \le 1  \},
\enas 
and for random vectors $\mathbf{X}$ and $\mathbf{Y}$ in $\mathbb{R}^p$, 
define the smooth functions metric
\bea \label{smoothdef}
d_{\mathcal{H}_{m,\infty,p}}\big(\mathcal{L}(\mathbf{X}),\mathcal{L}(\mathbf{Y})\big) = \sup_{h \in {{\cal H}}_{m,\infty,p}} |\E h(\mathbf{X})-\E h(\mathbf{Y})|.
\ena

For a positive semidefinite matrix $H$ we let $H^{1/2}$ denote the unique positive semidefinite square root of $H$. When $H$ is positive definite, we write $H^{-1/2}=(H^{1/2})^{-1}$. Our main multidimensional result is the following theorem:

\begin{theorem} \label{stein:fkg2}
With $m,p \in \mathbb{N}_1$, let $\{ \xi_{i,j}:i \in [m],j\in [p] \}$ be positively associated mean zero random variables bounded in absolute value by some positive constant $B$. Let $\mathbf{S} = (S_1,S_2,\ldots,S_p)$ where $S_j= \sum_{1 \le i \le m}\xi_{i,j}$ for $j \in [p]$ and assume that $\Sigma = \Var \big({\bf S} \big)$ is positive definite. Then
\begin{multline} \label{steinbound:fkg2}
d_{\mathcal{H}_{3,\infty,p}}\big(\mathcal{L}(\Sigma^{-1/2} \mathbf{S}),\mathcal{L}(\mathbf{Z})\big)\le 
\left(\frac{1}{6}+2\sqrt{2}\right)p^3B|\Sigma^{-1/2}|_{\infty}^3\sum_{j=1}^p\Sigma_{j,j} \\+ \left(\frac{3}{\sqrt{2}}+\frac{1}{2}\right)p^2|\Sigma^{-1/2}|_{\infty}^2 \sum_{j=1}^p \sum_{i,k\in[m],i \ne k} \Cov\left(\xi_{i,j},\xi_{k,j}\right) \\
+ \left( 2\sqrt{2}p^3B |\Sigma^{-1/2}|_{\infty}^3 + \left(\frac{3}{\sqrt{2}}+\frac{1}{2}\right)p^2|\Sigma^{-1/2}|_{\infty}^2 \right) \sum_{j,l \in [p],j \ne l} \Sigma_{j,l},
\end{multline}
where $\mathbf{Z} \sim {\cal N}({\bf 0},I_p)$, a standard normal vector in $\mathbb{R}^p$.
\end{theorem}

The Ising and percolation models are special cases of random cluster models, see \cite{Gri06}, for which positive association is well-known. These two cases, however, are the only ones for which covariances are known to decay exponentially. However, as in some generality the boundedness hypothesis of our main theorems hold for random cluster models, our results apply there whenever covariances can be shown to be sufficiently small. In particular, under known exponential decay, Corollary \ref{cor:isingpercolation} may be applied; see Remark \ref{re:AD15.etc}.

On the other hand, our results for the voter model in Section \ref{app3} illustrate that our methods apply in the absence of exponential decay. For Potts models, corresponding to the parameters $q \in \mathbb{N}_2$ of the random cluster model, it would therefore be of interest to determine if the covariance formula given in Theorem 9.11 of \cite{Gri06} admits sufficiently rapid covariance decay that would allow our methods to apply.  
We note that neither the voter model nor the contact process studied in Sections \ref{app3} and \ref{app4}, respectively, are instances of the random cluster model; we apply specific techniques to handle their covariance terms.

The remainder of this work is organized as follows. In Section \ref{app1} we state general results for positively associated random fields whose covariance decays exponentially. These results are then applied to the Ising and percolation models in Section \ref{app2}. The two main general results of this section, Theorems \ref{stein:fkg} and \ref{stein:fkg2}, are applied directly to the voter model and the contact process in Sections \ref{app3} and \ref{app4}. We then prove the general results of this section in Section \ref{pf}. One advantage of these latter results are that, unlike many results based on Stein's method, they may be applied without the need for coupling constructions.

\section{Applications} %\label{app}
We begin this section with a general result for a certain class of positively associated random fields.

\subsection{Second order stationary positively associated fields with exponential covariance decay} \label{app1}
Let $\{ X_{{\bf j}}: {\bf j}\in \mathbb{Z}^d \}$ be a positively associated random field on the $d$-dimensional integer lattice $\mathbb{Z}^d$. Assume that the field is second order stationary, that is, for all ${\bf j},{\bf k} \in \mathbb{Z}^d$ the covariance ${\rm Cov}(X_{{\bf j}},X_{{\bf k}})$ exists and equals 
$R({\bf j}-{\bf k})$ for some function $R$. With $S_{\bf k}^n$ defined in \eqref{blockdef}, consider the standardized variables
\bea \label{wn2}
W_{{\bf k}}^n = \frac{S_{{\bf k}}^n - \E S_{{\bf k}}^n}{\sqrt{n^dA_n}}, \qmq{${\bf k} \in \mathbb{Z}^d, n \in \mathbb{N}_1$,}
\ena
that have mean zero and variance 1. The following theorem provides a bound of order $n^{-d/(2d+2)}$ with an explicit constant
on the $L^1$ distance between the distribution of $W_{{\bf k}}^n$ and the normal under the condition that the covariance function $R(\cdot)$ decays exponentially in the $L^1$ norm in $\mathbb{R}^d$. We have chosen to bound the covariance decay in the $L^1$ norm for convenience in the proof of Lemma \ref{lemmaising}; as all norms in $\mathbb{R}^d$ are equivalent, when inequality \eqref{Ck.exp.bound2} holds for any norm with some $\lambda>0$ then it holds for any other norm with $\lambda$ replaced by some positive constant multiple.

\begin{theorem} \label{multiising} 
Let $d \in \mathbb{N}_1$ and  $\{ X_{{\bf j}}: {\bf j}\in \mathbb{Z}^d \}$ be a positively associated second order stationary random field with covariance function $R({\bf k})={\rm Cov}(X_{{\bf j}},X_{{\bf j}+{\bf k}})$ for all ${\bf j},{\bf k} \in {\mathbb Z}^d$, and 
suppose that for some $K > 0$ it holds that $|X_{{\bf j}}| \le K$ a.s.\ for all ${\bf j} \in \mathbb{Z}^d$. Assume that there exist $\lambda >0$ and $\kappa_0 >0$ such that 
\bea \label{Ck.exp.bound2}
R({\bf k}) \le \kappa_0 e^{-\lambda |{\bf k}|_1} \qmq{for all ${\bf k} \in \mathbb{Z}^d$ }
\ena
and let
\bea \label{muandupsilon}
\mu_{\lambda} = \frac{e^{\lambda}}{\left(e^{\lambda}-1\right)^2} \text{, \ } \upsilon_{\lambda}= \frac{e^{2\lambda}}{\left(e^{\lambda}-1\right)^2} \qmq{and} \gamma_{\lambda,d}= (4\mu_\lambda+2\upsilon_{\lambda})^d-(2\upsilon_{\lambda})^d
\ena
and	
\beas
C_{\lambda,\kappa_0,d}=\frac{5Kd\sqrt{\pi A_n}}{\sqrt{2}\kappa_0\gamma_{\lambda,d}}.
\enas 
Then, for any ${\bf k} \in \mathbb{Z}^d$, with $W_{{\bf k}}^n$ as given in \eqref{wn2} and $Z$ a standard normal random variable, 
\beas
d_1\big({\cal L}(W_{{\bf k}}^n),{\cal L}(Z)\big) \le \frac{\kappa_1}{n^{d/(2d+2)}} \text{ \ for all \ } n \ge \max \left\{ C_{\lambda,\kappa_0,d}^{2/d},C_{\lambda,\kappa_0,d}^{-2/(d+2)} \right\},
\enas
where, with $A_n$ is as in \eqref{andef},
\bea \label{kappavalue2} 
\kappa_1 =  \left(\frac{10K\kappa_0^d
	\gamma_{\lambda,d}^d 2^{3d/2}}{\pi^{d/2}A_n^{d+1/2}}\right)^{1/(d+1)} \left(\frac{1
}{d^{\frac{d}{d+1}}}+ 2d^{\frac{1}{d+1}}\right). 
\ena

\end{theorem}

\bigskip

In the following we will make use of the identities
\bea \label{Sum.krk4}
\sum_{k=1}^{n-1} (n-k)w^{k}= \frac{w\left((n-1)-nw+w^n\right)}{(w-1)^2} \text{ \ \ for \ \ } w \ne 1,
\ena
\bea \label{Sum.krk3}
n + \sum_{a=1}^{n-1} (n-a)(v^a+v^{-a})= \frac{v^{1-n} \left(v^n-1\right)^2}{(v-1)^2} \text{ \ \ for \ \ } v \ne 1,
\ena
and
\bea \label{Sum.krk2}
n+ 2\sum_{b=1}^{n-1} (n-b)u^{b}= \frac{(1-u^2)n-2u+2u^{n+1}}{(u-1)^2} \text{ \ \ for \ \ } u \ne 1.
\ena

%A<\infty was formerly \eqref{newmancondition}
We have already taken $A$ in \eqref{def:A}, without loss of generality, to be positive, and the exponential decay condition \eqref{Ck.exp.bound2} implies that it is also finite. Lemma 4 of \cite{new80} can now be invoked to yield that $\lim_{n \rightarrow \infty} A_n = A$. Hence the values $A_n$ are bounded away from zero and infinity as a function of $n$, and replacing $A_n$ by its limiting value $A$ in \eqref{kappavalue2} only affects the bound by a constant and does not alter the rate of convergence. 

In certain instances the values of $A_n$ and $A$ will be close even for moderate the values $n$, for instance if equality holds in \eqref{Ck.exp.bound2}, then using \eqref{Sum.krk2}
with $u = e^{-\lambda}$, which is not equal to one since $\lambda > 0$, as $n \rightarrow \infty$ we obtain
\begin{multline*}
A_n = \frac{\kappa_0}{n^d}  \sum_{{\bf i},{\bf j} \in B_{{\bf 1}}^n} e^{-\lambda|{\bf i}-{\bf j}|_1} 
    = \frac{\kappa_0}{n^d} \sum_{\substack{i_1,\ldots,i_d =1 \\ j_1,\ldots,j_d=1}}^n \prod_{k=1}^d e^{-\lambda|i_k-i_k|} 
		= \frac{\kappa_0}{n^d} \prod_{k=1}^d \sum_{i_k,j_k=1}^n e^{-\lambda|i_k-i_k|} \\
		= \frac{\kappa_0}{n^d} \prod_{k=1}^d \sum_{a_k=-n+1}^{n-1} (n-|a_k|) e^{-\lambda |a_k|} 
		= \frac{\kappa_0}{n^d} \left(n+2\sum_{b=1}^{n-1} (n-b)e^{-\lambda b}\right)^{d} \\
		= \kappa_0\left( \frac{1-e^{-2\lambda} -\frac{2e^{-\lambda}}{n}+\frac{2e^{-\lambda (n+1)}}{n}}{\left(1-e^{-\lambda}\right)^2}  \right)^d \longrightarrow \kappa_0\left( \frac{1+e^{-\lambda}}{1-e^{-\lambda}}  \right)^d
		=\kappa_0 \coth^d(\lambda/2):=A. %\quad \mbox{as $n \rightarrow \infty$.}
\end{multline*}

One can apply Theorem \ref{Bultheorem} above, due to \cite{Bul95}, to obtain a bound of order $\log^d n/\sqrt{n}$, without an explicit constant, on the Kolmogorov distance between $W_{{\bf k}}^n$ and the normal
by showing that $u(n)$ as defined in \eqref{undef} converges to zero at an exponential rate. Here we use Theorem \ref{stein:fkg} to obtain an $L^1$ result.

We note that the technique in \cite{Bul95}, first introduced in \cite{Tik80}, presents some challenges when applied to obtain an $L^1$ bound. We consider the argument in the one dimensional case, as in \cite{Bir88}, to illustrate the difficulty. Letting $f_n$ be the characteristic function of $W_1^n$ and following the proof in \cite{Bir88}, one obtains the bound
\beas
\left|f_n(t)-e^{-t^2/2}\right| \le C\left(n^{-1/2}\log (n) t^3 e^{-t^2/4} + n^{-1}\log^{3/2}(n) t + n^{-3/2}t\right)
\enas
over the region $0 \le t \le \gamma n^{1/2}/\log (n)$, where $\gamma$ is a constant defined in \cite{Bir88}, and where here and in the following, $C$ denotes a constant whose value may change from line to line. Applying Theorem 1.5.2 in \cite{IL71} with $T=n^{1/2}/\log (n) $, the Kolmogorov distance of order $\log (n)/n^{1/2}$ was obtained in \cite{Bir88}. To get the bound in $L^1$, one may apply Theorem 1.5.4 in \cite{IL71} instead, however the integrand of the term (3) in Theorem 1.5.4 
\beas
\int_{-T}^T \left(\frac{d}{dt}\frac{f_n(t)-e^{-t^2/2}}{t} \right)^2 dt \quad \mbox{includes} \quad \frac{\left(f_n(t)-e^{-t^2/2}\right)^2}{t^4},
\enas 
which is not integrable.

We also extend Theorem \ref{multiising} to the multidimensional case. For any $p \in \mathbb{N}_1$ and indices ${\bf k}_1,\ldots,{\bf k}_p$ in $\mathbb{Z}^d$ satisfying a separation condition to ensure non-degeneracy of the limiting distibution, Theorem \ref{multiising2} provides a bound in the metric $d_{\mathcal{H}_{3,\infty,p}}$ to the multivariate normal for ${\bf S}^n=(S_{{\bf k}_1}^n,\ldots,S_{{\bf k}_p}^n)$ under exponential decay of the covariance function and 
a condition limiting the amount of variables common to the coordinates making up the sums ${\bf S}^n$. 
Regarding multidimensional convergence, under strictly stationary Theorem \ref{newmantheorem} of \cite{new80} shows that $(W_{n{\bf k}_1}^n,\ldots,W_{n{\bf k}_p}^n)$ converges to a standard normal random vector in $\mathbb{R}^p$, but does not provide a bound. In the following result and its proof, constants will not be tracked with precision, but will be indexed by the set of variables on which it depends.

\begin{theorem} \label{multiising2}
For $d \in \mathbb{N}_1$, let $\{ X_{{\bf j}}: {\bf j}\in \mathbb{Z}^d \}$ be a positively associated second order stationary random field with covariance function $R({\bf k})={\rm Cov}(X_{{\bf j}},X_{{\bf j}+{\bf k}})$ for all ${\bf j},{\bf k} \in {\mathbb Z}^d$, and suppose that there exists constants $K>0,\kappa_0>0$ and $\lambda > 0$ such that $|X_{{\bf j}}| \le K$ a.s.\ for all ${\bf j} \in \mathbb{Z}^d$ and
\beas
R({\bf k}) \le \kappa_0e^{-\lambda |{\bf k}|_1} \text{ \ for all \ } {\bf k}\in \mathbb{Z}^d.
\enas
For $p\in \mathbb{N}_1$ let ${\bf k}_1,\ldots,{\bf k}_p \in \mathbb{Z}^d$ be 
such that
\bea \label{eq:kqks.separated.alphan}
\min_{q,s\in[p], q \ne s}\left|{\bf k}_q-{\bf k}_s\right|_{\infty} \ge (1-\alpha)n,
\ena
for some $0<\alpha<1$ with $\alpha n$ an integer. Let
$\mathbf{S}^n = (S_{{\bf k}_1}^n ,\ldots,S_{{\bf k}_p}^n )$, where $S_{{\bf k}}^n$ is defined as in \eqref{blockdef} and assume that the covariance matrix $\Sigma$ of $\mathbf{S}^n$ is invertible, and let
\beas %\label{eq:defpsin}
\psi_n=n^{d/2}|\Sigma^{-1/2}|_\infty.
\enas 
Then, there exists a constant $C_{\lambda,\kappa_0,d,p,K}$ such that with $\mathbf{Z}$ a standard normal random vector in $\mathbb{R}^p$,
\begin{multline} \label{eq:thm2rangen}
d_{\mathcal{H}_{3,\infty,p}}\big(\mathcal{L}(\Sigma^{-1/2}(\mathbf{S}^n-\E \mathbf{S}^n)),\mathcal{L}(\mathbf{Z})\big) \\
\le C_{\lambda,\kappa_0,d,p,K} \left((A_n+\alpha)^{\frac{1}{d+1}} \psi_n^{\frac{2d+3}{d+1}}\left(\frac{1}{d^{\frac{d}{d+1}}}+2 d^{\frac{1}{d+1}}\right)n^{-\frac{d}{2(d+1)}}+\alpha \psi_n^2\right)\\
\mbox{for all} \quad n \ge \max \left\{B_{n,d,\alpha}^{2/d},B_{n,d,\alpha}^{-2/(d+2)}\right\}  \qmq{where} B_{n,d,\alpha}= d\psi_n (A_n+\alpha).
\end{multline}
\end{theorem}

For checking the hypothesis that the covariance matrix of ${\bf S}^n$ is invertible, and that the quantity $\psi_n$ is of order one, and so does not affect the rate in \eqref{eq:thm2rangen} and also so that $B_{n,d,\alpha}$ is order one, we present the following condition, proved at the end of this section. We note in addition that if $\psi_n=O(1)$ and $\alpha=O(n^{-\frac{d}{2(d+1)}})$ then the bound in \eqref{eq:thm2rangen} will also be of order $O(n^{-\frac{d}{2(d+1)}})$. 

\begin{lemma}\label{lem:Gershgorin}
Let $\{ X_{{\bf j}}: {\bf j}\in \mathbb{Z}^d \}$ be a second order stationary random field with covariance function $R({\bf k})={\rm Cov}(X_{{\bf j}},X_{{\bf j}+{\bf k}})$ for all ${\bf j},{\bf k} \in {\mathbb Z}^d$ where $R(\cdot)$ satisfies \eqref{Ck.exp.bound2}. With $n \in \mathbb{N}_2$, the covariance matrix $\Sigma$ of the random vector ${\bf S}^n \in \mathbb{R}^p, p \ge 2$ as given in Theorem \ref{multiising2} is invertible if for some $b>0$ 
\bea \label{eq:cond.b.invert}
\left|{\bf k}_q-{\bf k}_s\right|_{\infty} \ge n-b \qmq{for all $q \not =s, q,s \in [p]$, and} b< \frac{nA_n}{(p-1)\kappa_0\upsilon_{\lambda}^d},
\ena
with $\lambda$ and $\kappa_0$ as in \eqref{Ck.exp.bound2}, and $\upsilon_{\lambda}$ as in \eqref{muandupsilon}. If $b=\alpha n$ for $\alpha$ satisfying 
\bea \label{alphacondition}
0 < \alpha < \min\left\{1, \frac{A_n}{(p-1)\kappa_0\upsilon_{\lambda}^d} \right\}
\ena
then the matrix $\Sigma$ is invertible and
\bea \label{eq:infinity.norm.bound.inverse.lem:Gershgorin}
|\Sigma^{-1}|_{\infty} \le \frac{1}{n^d(A_n - (p-1)\kappa_0 \upsilon_{\lambda}^d\alpha)}.
\ena
\end{lemma}

The proofs of Theorems \ref{multiising}  and \ref{multiising2} proceed by decomposing the sum $S_{\bf k}^n$ over the block $B_{\bf k}^n$
into sums over smaller, disjoint blocks whose side lengths are at most some integer $l$. In particular, 
given $1 \le l \le n$ uniquely write $n=(m-1)l+r$ for $m \ge 1$  and $1 \le r \le l$. We correspondingly decompose $B_{{\bf k}}^n$ into $m^d$ disjoint blocks $D_{{\bf i},{\bf k}}^l, {\bf i} \in [m]^d$, where there are $(m-1)^d$ `main' blocks having all sides of length $l$, and $m^d-(m-1)^d$ remainder blocks having all sides of length $r$ or $l$, with at least one side of length $r$.

In detail, for ${\bf k} \in \mathbb{Z}^d$ and ${\bf i} \in [m]^d$ set $D_{{\bf i},{\bf k}}^l=D_{\bf i}^l+{\bf k} - {\bf 1}$ where
\beas
D_{\bf i}^l=\{{\bf j} \in \mathbb{Z}^d: (i_s-1)l+1 \le j_s \le i_sl, i_s \not =m,  (m-1)l+1 \le j_s \le (m-1)l+r, i_s =m\}.
\enas
It is straightforward to verify that for all ${\bf i} \in [m]^d$ with no component equal to $m$, which are the indices of the `main blocks', that
with $B_{{\bf 1}}^l$ as in \eqref{blockdef}, we have
\bea \label{Di.block}
D_{\bf i}^l = B_{{\bf 1}}^l+ ({\bf i}-{\bf 1})l \text{ \ \ for \ \ } {\bf i} \in [m-1]^d,
\ena
and if $r=l$ then $D_{{\bf i}}^l$ is given by \eqref{Di.block} for all $i \in [m]^d$. Further, it is easy to see that the elements of the collection $\{D_{{\bf i},{\bf k}}^l, {\bf i} \in [m]^d\}$ are disjoint, and union to $B_{\bf k}^n$.

Letting
\bea \label{xidef}
\xi_{{\bf i},{\bf k}}^{l} = \sum_{{\bf t} \in D_{{\bf i},{\bf k}}^l}\left(X_{{\bf t}} - \E X_{{\bf t}}\right) \quad \mbox{and} \quad W_{{\bf k}}^n = \sum_{{\bf i} \in [m]^d} \frac{\xi_{
		{\bf i},{\bf k}
		}^{l}}{\sqrt{n^d A_n}} \text{ \ \ \ for \ \ \  } {\bf i} \in [m]^d,
\ena 
we see that $\xi_{{\bf i},{\bf k}}^l$ has mean zero, and
$W_{{\bf k}}^n$ as in \eqref{xidef} agrees with its representation as given in \eqref{wn2}, and has mean zero and variance one. For simplicity we will drop the index ${\bf k}$ in $\xi_{{\bf i},{\bf k}}^l$ when ${\bf k}={\bf 1}$, as we do also for $D_{{\bf i},{\bf k}}$, and may also suppress $l$. 

By \cite{EPW67}, when ${\bf X}=(X_1,\ldots,X_p)$ is positively associated and if $h_k:\mathbb{R}^p \rightarrow \mathbb{R}$, $k \in [q]$ are coordinate-wise nondecreasing functions, then $h_1({\bf X}),\ldots, h_q({\bf X})$ are positively associated. Hence, as the elements of $\{ \xi_{{\bf i},{\bf k}}: {\bf i} \in [m]^d \}$ are increasing functions of $\{ X_{{\bf j}}: {\bf j} \in \mathbb{Z}^d \}$ they are positively associated. 

We prove Theorems \ref{multiising} and \ref{multiising2} with the help of two lemmas that follow.  The first, Lemma \ref{lemmaising} bounds the sum of the covariances between $\xi_{{\bf i},{\bf k}}^{l}$ and $\xi_{{\bf j},{\bf k}}^{l}$, defined in \eqref{xidef}, over ${\bf i},{\bf j} \in [m]^d$. Next, Lemma \ref{lemmaising2} bounds the covariance between two block sums of size $n^d$ whose centers are at least $n-b$ apart in the $L_\infty$ distance for some integer $1 \le b \le n$. 

\begin{lemma} \label{lemmaising}
Let $\{ X_{{\bf j}}: {\bf j}\in \mathbb{Z}^d \}$ be a second order, positively associated stationary random field with covariance function $R({\bf k})={\rm Cov}(X_{{\bf j}},X_{{\bf j}+{\bf k}})$ for all ${\bf j},{\bf k} \in {\mathbb Z}^d$ where $R(\cdot)$ satisfies the exponential decay condition \eqref{Ck.exp.bound2}. For $n \ge 2$ and $1 \le l \le n$, let $n=(m-1)l+r$ for integers $m \in \mathbb{N}_1$ and $1 \le r \le l$. Then for ${\bf k} \in \mathbb{Z}^d$, with $\xi_{{\bf i},{\bf k}}^l$ given by \eqref{xidef} we have 
\beas
\sum_{{\bf i},{\bf j} \in [m]^d,{\bf i} \neq {\bf j}} \E \left[\xi_{{\bf i},{\bf k}}^{l} \xi_{{\bf j},{\bf k}}^{l}\right] \le \frac{\kappa_0\gamma_{\lambda,d}n^d}{l},
\enas
where $\kappa_0$ is given in \eqref{Ck.exp.bound2}, and $\gamma_{\lambda,d}$ in \eqref{muandupsilon}. 
\end{lemma}

\proof By second order stationarity, it suffices to consider the case ${\bf k}={\bf 1}$. For convenience write
$\sigma_{{\bf i},{\bf j}}=\E \left[\xi_{{\bf i}} \xi_{{\bf j}}\right]$.
Since $D_{{\bf i}}^l \subset B_{{\bf 1}}^l+ ({\bf i}-{\bf 1})l$ for all ${\bf i} \in [m]^d$ and $\Cov\big(X_{{\bf t}},X_{{\bf s}}\big) \ge 0$ for all ${\bf t}, {\bf s} \in \mathbb{Z}^d$ by positive association, with $S_{{\bf i}}^l$ as given in \eqref{blockdef}, we have
\bea \label{remainderub}
\sigma_{{\bf i},{\bf j}} = \Cov \left(\sum_{{\bf t} \in D_{{\bf i}}^l}X_{{\bf t}} ,\sum_{{\bf s} \in D_{{\bf j}}^l}X_{{\bf s}} \right)\le \Cov \big(S_{{\bf 1}+({\bf i}-{\bf 1})l}^l ,S_{{\bf 1}+({\bf j}-{\bf 1})l}^l \big) \text{ \ \ for all \ \ } {\bf i} \ne {\bf j} \in [m]^d.
\ena
Note that when ${\bf i} \ne {\bf j}$ in $[m]^d$ the pair $({\bf i},{\bf j})$ must lie in exactly one of the sets $E_s^m$ for $s\in [d]$, given by
\beas %\label{esdef}
E_s^m=\left\{({\bf i},{\bf j}) \in [m]^d \times [m]^d: |\{k \in [d]: i_k \not = j_k \} |=s \right\} \quad \mbox{for $s\in [d]$.}
\enas

Hence, using \eqref{remainderub},
\bea \label{lemmaclaim1}
\sum_{{\bf i},{\bf j} \in [m]^d,{\bf i} \neq {\bf j}}\sigma_{{\bf i},{\bf j}} &=& \sum_{s=1}^d \sum_{({\bf i},{\bf j}) \in E_s^m  } \sigma_{{\bf i},{\bf j}} \le \sum_{s=1}^d \sum_{({\bf i},{\bf j}) \in E_s^m  } \Cov \big(S_{{\bf 1}+({\bf i}-{\bf 1})l}^l ,S_{{\bf 1}+({\bf j}-{\bf 1})l}^l \big).
\ena

Recalling the definition of the block sums in \eqref{blockdef} and using second order stationarity, for $({\bf i},{\bf j}) \in E_s^m$ we have

\begin{multline}
\Cov \left(S_{{\bf 1}+({\bf i}-{\bf 1})l}^l,S_{{\bf 1}+({\bf j}-{\bf 1})l}^l \right) 
= \sum_{\substack{p_1,\ldots,p_d =1\\  q_1,\ldots,q_d=1}}^{l} R\left(\begin{bmatrix} (p_1+(j_1-1)l)-(q_1+(i_1-1)l) \\ \vdots \\ (p_d+(j_d-1)l)-(q_d+(i_d-1)l) \end{bmatrix}\right)  \\ 
= \sum_{\substack{p_1,\ldots,p_d =1\\ q_1,\ldots,q_d=1}}^{l} R\left(\begin{bmatrix} p_1-q_1+(j_1-i_1)l \\ \vdots \\ p_d-q_d+(j_d-i_d)l \end{bmatrix}\right)  \\ 
%&=&   \sum_{a_1,\ldots,a_d=-l+1}^{l-1} (l-|a_1|)\ldots(l-|a_d|)R\left(\begin{bmatrix} (j_1-i_1)l+a_1\\ \vdots \\ (j_d-i_d)l+a_d \end{bmatrix}\right) \\
\le \kappa_0 \sum_{a_1,\ldots,a_d=-l+1}^{l-1} (l-|a_1|)\cdots(l-|a_d|)\exp\left(-\lambda\left|\begin{bmatrix} a_1+(j_1-i_1)l\\  \vdots \\ a_d+(j_d-i_d)l  \end{bmatrix}\right|_1\right), \label{eq:l1.before.pi}
\end{multline}
where we have applied \eqref{Ck.exp.bound2} in the final inequality.

We now apply some invariance properties of the norm $|{\bf x}|_1$ in order to simplify the sum of expression \eqref{eq:l1.before.pi} when taken over pairs of indices in $E_s^m$. We say that the pairs of indices $({\bf i},{\bf j})$ and $({\bf i}',{\bf j}')$ in $E_s^m$ are equivalent if $({\bf i},{\bf j})$ can be made to agree with $({\bf i}',{\bf j}')$ by interchanging $i_k$ and $j_k$ for any values of $k \in [d]$ and then applying the same permutation to the resulting set of vectors.  Note that as $a_k$ and $-a_k$ appear symmetrically in the sum \eqref{eq:l1.before.pi}, and that $|\cdot|_1$ is invariant under sign changes and permutations of coordinates, equivalent vectors yield the same value of this sum.

Call an index pair $({\bf i},{\bf j}) \in E_s^m$ canonical if $i_k< j_k$ for all $k=1,\ldots,s$. Note that as pairs in $E_s^m$ have exactly $s$ inequalities among their coordinates that canonical indices must agree in their remaining, and last, $d-s$ coordinates. Now associate to every pair $({\bf i},{\bf j}) \in E_s^m$ the equivalent, unique canonical vector that is obtained from $({\bf i},{\bf j})$ by interchanging all coordinates $k \in [d]$ for which $i_k>j_k$, and then applying the permutation that maps the $s$ unequal coordinates of the pair to positions $1,2,\ldots,s$, thus mapping the  equal pairs to the last $d-s$ positions, leaving the relative order of each group unchanged. Partitioning the sum over $E_s^m$ into smaller sums over all vectors that are equivalent to the same canonical vector yields a sum over all canonical vectors, with a factor of $2^s$ to account for the number arrangements of inequalities over the unequal pairs, and a factor of ${d \choose s}$ to account for the positions in which the unequal coordinates may occur. Hence, summing \eqref{eq:l1.before.pi} over $E_s^m$, and using identities of the form
\beas
\sum_{a_1\in A_1,\ldots,a_d \in A_d}\prod_{q=1}^d f_q(a_q)=\prod_{q=1}^d \sum_{a_q\in A_q} f_q(a_q)
\enas
that hold for any index sets $I_1,\ldots,I_d$ and any functions $f_q$ in order to interchange summations and products notating the index set to be consistent with the formula used above
\beas
I_d^{s,m} = \{({\bf i},{\bf j}) \in \mathbb{Z}^d\times \mathbb{Z}^d: 1\le i_k < j_k \le m,k=1,\ldots,s,1\le i_k = j_k \le m,k=s+1,\ldots,d \},
\enas
we obtain

\beas
&&\sum_{({\bf i},{\bf j}) \in E_s^m  }\Cov \left(S_{{\bf 1}+({\bf i}-{\bf 1})l}^l,S_{{\bf 1}+({\bf j}-{\bf 1})l}^l \right) \\
&\le&  \kappa_0 {d \choose s} 2^s \sum_{({\bf i},{\bf j}) \in I_d^{s,m}}  \sum_{a_1,\ldots,a_d=-l+1}^{l-1} (l-|a_1|)\cdots(l-|a_d|)\exp\left(-\lambda\left|\begin{bmatrix} a_1+(j_1-i_1)l\\  \vdots \\ a_d+(j_d-i_d)l \end{bmatrix}\right|_1\right) \\
&=& \kappa_0 {d \choose s} 2^s \sum_{({\bf i},{\bf j}) \in I_d^{s,m}}  
 \sum_{a_1,\ldots,a_d=-l+1}^{l-1} \prod_{q=1}^d (l-|a_q|) e^{-\lambda|(j_q-i_q)l+ a_q|}\\
&=& \kappa_0 {d \choose s} 2^s \sum_{({\bf i},{\bf j}) \in I_d^{s,m}} 
	\prod_{q=1}^d \sum_{a_q=-l+1}^{l-1} (l-|a_q|) e^{-\lambda|(j_q-i_q)l+ a_q|}\\
&=& \kappa_0 {d \choose s} 2^s  \sum_{({\bf i},{\bf j}) \in I_d^{s,m}} 
	\left(\prod_{q=1}^s \sum_{a_q=-l+1}^{l-1} (l-|a_q|) e^{-\lambda(j_q-i_q)l-\lambda a_q}\right)\left(\prod_{q=s+1}^d \sum_{a_q=-l+1}^{l-1} (l-|a_q|)e^{-\lambda |a_q|}\right),
\enas
where in the final equality we have expressed the product in $q$ as two separate products, and have noted in the first of these, as $j_k-i_k \ge 1$ and $a_k$ ranges from $-l+1$ to $l-1$, the terms $(j_k-i_k)l+a_k$ appearing in the absolute value of the exponent are always positive. Now expanding the first multiple summation as an inner summation on the first $s$ indices and outer sum over the final $d-s$ indices, separating out the exponential terms in the first product and noting that the remaining product terms are purely mulitplicative factors, we obtain

\beas
&& \kappa_0 {d \choose s} 2^s \sum_{\substack{1 \le i_k = j_k \le m \\ k=s+1,\ldots,d}}   
                          \sum_{\substack{1 \le i_k < j_k \le m \\ k=1,\ldots,s}} \prod_{q=1}^s e^{-\lambda(j_q-i_q)l}\\
          && \text{\hspace{100pt}} \times\left(\prod_{q=1}^s \sum_{a_q=-l+1}^{l-1} (l-|a_q|) e^{-\lambda a_q}\right) 
					                                                   \left(\prod_{q=s+1}^d \sum_{a_q=-l+1}^{l-1} (l-|a_q|)e^{-\lambda |a_q|} \right)\\
&=& \kappa_0 {d \choose s} 2^s \left(\sum_{i=1}^m 1\right)^{d-s}  \left(\prod_{q=1}^s\sum_{1 \le i_q < j_q \le m}  e^{-\lambda(j_q-i_q)l}\right) \\
          && \text{\hspace{100pt}} \times\left(\prod_{q=1}^s \sum_{a_q=-l+1}^{l-1} (l-|a_q|) e^{-\lambda a_q}\right) 
                                                            \left(\prod_{q=s+1}^d \sum_{a_q=-l+1}^{l-1} (l-|a_q|)e^{-\lambda |a_q|} \right)\\
&=& \kappa_0 {d \choose s} 2^s m^{d-s}  \left(\sum_{k=1}^{m-1}(m-k)e^{-\lambda k l}\right)^s\\
          && \text{\hspace{100pt}} \times\left(l+\sum_{a=1}^{l-1} (l-a) \left(e^{\lambda a}+e^{-\lambda a}\right)\right)^s\left(l+2\sum_{b=1}^{l-1} (l-b)e^{-\lambda b}\right)^{d-s}.
\enas 

By applying identity \eqref{Sum.krk4} with $n=m$ and $w = e^{-\lambda l}$ and identities \eqref{Sum.krk3} and \eqref{Sum.krk2} with $n=l$, $v = e^{\lambda}$ and $u = e^{-\lambda}$, noting that $u$, $v$ and $w$ are not one as $\lambda > 0$, the term above equals

\bea \label{lemmaclaim2}
&&\kappa_0 {d \choose s} 2^s m^{d-s} \left( \frac{e^{\lambda}e^{-2\lambda l} \left(e^{\lambda l}-1\right)^2\left(m-1 -me^{-\lambda l}+e^{-\lambda l m}\right)}{(e^{\lambda}-1)^2(1-e^{-\lambda l})^2} \right)^s \nn \\
          && \text{\hspace{150pt}} \times 
					\left( \frac{\left(1-e^{-2\lambda}\right)l -2e^{-\lambda}+2e^{-\lambda (1+l)}}{\left(1-e^{-\lambda}\right)^2}  \right)^{d-s}  \nn \\
&\le& \kappa_0 {d \choose s} 2^s m^d  \left(\frac{e^{\lambda}}{\left(e^{\lambda}-1\right)^2}\right)^s\left(\frac{l}{\left(1-e^{-\lambda}\right)^2}\right)^{d-s} \nn\\
&=& \kappa_0 m^d {d \choose s} (2\mu_{\lambda})^s \left(l\upsilon_{\lambda}\right)^{d-s} 
\ena 
where $\mu_{\lambda}$ and $\upsilon_{\lambda}$ are given in \eqref{muandupsilon}. By \eqref{lemmaclaim1} and \eqref{lemmaclaim2}, we have
\beas %\label{boundis3}
\sum_{{\bf i},{\bf j} \in [m]^d,{\bf i} \neq {\bf j}} \sigma_{{\bf i},{\bf j}} 
&\le& \kappa_0m^d \sum_{s=1}^d {d \choose s}  (2\mu_{\lambda})^s\left(l\upsilon_{\lambda}\right)^{d-s} \nn \\
&\le& \frac{\kappa_0 n^d}{l}\sum_{s=1}^d {d \choose s} (4\mu_{\lambda})^s\left(2\upsilon_{\lambda}\right)^{d-s} \nn \\
&=&\frac{\kappa_0n^d\left((4\mu_{\lambda}+2\upsilon_{\lambda})^d-\left(2\upsilon_{\lambda}\right)^d\right)}{l} 
\enas
where we have used the bounds $m \le 2n/l$ and $1/l^s \le 1/l$ for all $s \in [d]$ in the second inequality. \bbox

\bigskip
\begin{lemma} \label{sumexpo}
For all $n \in \mathbb{N}_2$ and $\lambda > 0$,
\beas
\sum_{a=-n+1}^{n-1} (n-|a|)e^{-\lambda |q+a|}  \text{ \ is decreasing as a function of $|q| \in \mathbb{N}_0$. }
\enas
\end{lemma}
\proof We note that the sequence $f(a)=(n-|a|){\bf 1}(|a| \le n)$ is unimodal, and $g(a)=e^{-\lambda|a|}$ is log-concave, and hence their convolution 
\beas
h(q)=\sum_{a=-\infty}^\infty f(a)g(q-a)=\sum_{a=-n+1}^{n-1} (n-|a|)e^{-\lambda|q-a|}=\sum_{a=-n+1}^{n-1} (n-|a|)e^{-\lambda|q+a|}
\enas
is unimodal, see \cite{Bre89}. Since $g(a)$ is symmetric unimodal, one can also use the fact from \cite{Sta89} that the convolution of two symmetric unimodal sequences is symmetric unimodal to prove the result. 
\bbox

\bigskip

Below we index the coordinates of vectors ${\bf k}_j \in \mathbb{Z}^d$ as $\left(k_1^j,\ldots,k_d^j\right)$.

\begin{lemma} \label{lemmaising2}
Let $\{ X_{{\bf j}}: {\bf j}\in \mathbb{Z}^d \}$ be a second order stationary random field with covariance function $R({\bf k})={\rm Cov}(X_{{\bf j}},X_{{\bf j}+{\bf k}})$ for all ${\bf j},{\bf k} \in {\mathbb Z}^d$ where $R(\cdot)$ satisfies \eqref{Ck.exp.bound2}. Let $n \in \mathbb{N}_2$ and let $b$ be an integer such that $1\le b \le n$. Then if ${\bf k}_1$ and ${\bf k}_2$ are vectors in $\mathbb{Z}^d$ such that
\beas
\left|{\bf k}_1-{\bf k}_2\right|_{\infty} \ge n-b,
\enas
then with $\lambda$ and $\kappa_0$ as in \eqref{Ck.exp.bound2}, and $\upsilon_{\lambda}$ as in \eqref{muandupsilon},
\beas
\Cov \left(S_{{\bf k}_1}^n,S_{{\bf k}_2}^n \right) \le \kappa_0 \upsilon_{\lambda}^d bn^{d-1}.
\enas 
\end{lemma}

\proof
Arguing as in the proof of Lemma \ref{lemmaising} we have
\bea
\Cov \left(S_{{\bf k}_1}^n,S_{{\bf k}_2}^n \right)&=&\sum_{\substack{p_1,\ldots,p_d =0 \\  q_1,\ldots,q_d=0}}^{n-1} R\left(\begin{bmatrix} (p_1+k_1^2)-(q_1+k_1^1) \\ \vdots \\ (p_d+k_d^2)-(q_d+k_d^1) \end{bmatrix}\right) \nonumber \\
&\le& \kappa_0 \sum_{a_1,\ldots,a_d=-n+1}^{n-1} (n-|a_1|)\ldots(n-|a_d|) \exp  \left(-\lambda\left|\begin{bmatrix} a_1+(k_1^2-k_1^1)\\ \vdots \\ a_d+(k_d^2-k_d^1) \end{bmatrix}\right|_1\right) \nonumber \\
&=& \kappa_0 \prod_{i=1}^d \sum_{a_i=-n+1}^{n-1} (n-|a_i|)e^{-\lambda |(k_i^2-k_i^1)+a_i|}.\label{eq:prod.term.one.exception}
\ena
Lemma \ref{sumexpo} yields that 
\bea \label{eq:dec.in.k1-k2}
\sum_{a_i=-n+1}^{n-1} (n-|a_i|)e^{-\lambda |(k_i^2-k_i^1)+a_i|} \qmq{is a decreasing function of $|k_i^1-k_i^2|$.}
\ena
In particular the $i^{th}$ sum appearing in the product
\eqref{eq:prod.term.one.exception}	
is bounded by its value when $k_i^1=k_i^2$. As $\left|{\bf k}_1-{\bf k}_2\right|_{\infty} \ge n-b$, there must exist at least one $i$ for which  $|k_i^2-k_i^1| \ge n-b$, and whose corresponding sum is maximized by its value when equality to $n-b$ is achieved, again using \eqref{eq:dec.in.k1-k2}. The product of these sums, by \eqref{eq:dec.in.k1-k2} again, is maximized when there is just a single coordinate achieving $n-b$ as its absolute difference, and where this difference in all other terms achieve equality to zero.  Hence, by symmetry the term above is bounded by the case where $k_i^1=k_i^2$ for $i \in [d-1]$ and $k_d^2-k_d^1=n-b$ and thus
\bea 
\Cov \left(S_{{\bf k}_1}^n,S_{{\bf k}_2}^n \right) 
&\le& \kappa_0\prod_{i=1}^{d-1} \sum_{a_i=-n+1}^{n-1} (n-|a_i|) e^{-\lambda|a_i|} \sum_{a_d=-n+1}^{n-1} (n-|a_d|)e^{-\lambda|a_d+n-b|} \nn\\
&\le& \kappa_0\left(n\upsilon_{\lambda} \right)^{d-1} \sum_{a_d=-n+1}^{n-1} (n-|a_d|)e^{-\lambda |a_d+n-b|},\label{sumlemma2.5}
\ena
where we have applied \eqref{Sum.krk2} in the final inequality and $\upsilon_{\lambda}$ is given in \eqref{muandupsilon}.

Now considering the sum in \eqref{sumlemma2.5} for $2 \le b \le n$, we obtain
\beas
\sum_{a=-n+1}^{n-1} (n-|a|)e^{-\lambda |a+n-b|} &=& \sum_{a=-n+1}^{-n+b-1} (n-|a|) e^{\lambda (a+n-b)} +\sum_{a=-n+b}^{n-1}(n-|a|) e^{-\lambda (a+n-b)} \\
&=& \sum_{a=n-b+1}^{n-1}(n-a)e^{\lambda (-a+n-b)} + \sum_{a=1}^{n-b} (n-a)e^{-\lambda (-a+n-b)}  \\
&& \text{\hspace{50pt}}+ \sum_{a=0}^{n-1}(n-a)e^{-\lambda (a+n-b)}.
\enas

For the first sum, making a change of variable and applying \eqref{Sum.krk4}, we obtain
\begin{multline*}
\sum_{a=n-b+1}^{n-1}(n-a)e^{\lambda (-a+n-b)} = e^{-\lambda b}\sum_{a=1}^{b-1}ae^{\lambda a}
= e^{-\lambda b}\left( b\sum_{a=1}^{b-1}e^{\lambda a}-\sum_{a=1}^{b-1}(b-a)e^{\lambda a}\right)\\ =e^{-\lambda b} \left(b\frac{e^{\lambda b}-e^\lambda}{e^\lambda-1} - \frac{e^\lambda((b-1)-be^\lambda+e^{\lambda b})}{(e^\lambda -1)^2}\right)\\
=\frac{b\left(e^{\lambda}-1\right)+e^{\lambda(1-b)}-e^{\lambda}}{\left(e^{\lambda}-1\right)^2},
\end{multline*}
and similarly we may obtain
\beas
\sum_{a=1}^{n-b} (n-a)e^{-\lambda (-a+n-b)} = \frac{b e^{\lambda}\left(e^{\lambda}-1\right)+e^{\lambda}+n\left(1-e^{\lambda}\right)e^{\lambda(1+b-n)}-e^{\lambda(1+b-n)}}{\left(e^{\lambda}-1\right)^2}
\enas
and
\beas
\sum_{a=0}^{n-1}(n-a)e^{-\lambda (a+n-b)} = \frac{e^{\lambda(1+b-2n)}-e^{\lambda(1+b-n)}+n(e^{\lambda}-1)e^{\lambda(1+b-n)}}{\left(e^{\lambda}-1\right)^2}.
\enas
Summing these three terms yields
\beas
\sum_{a=-n+1}^{n-1} (n-|a|)e^{-\lambda |a+n-b|}
=\frac{b (e^{2 \lambda }-1)+e^{(1-b) \lambda }-e^{-\lambda  (n-b-1)}(2-e^{-\lambda n})}{\left(e^{\lambda }-1\right)^2}
\le \frac{be^{2 \lambda}}{\left(e^{\lambda }-1\right)^2}
=b \upsilon_{\lambda}, 
\enas
where we lastly note that the equality holds also for $b=1$. 

\bbox 

\bigskip

Now we use Theorem \ref{stein:fkg} and Lemma \ref{lemmaising} to prove Theorem \ref{multiising}. In the following, for positive numbers $a$ and $b$ we will seek to minimize a quantity of the form $al^d +b/l$ over non-negative integers $l$. Over real values, it is easy to verify that the minimum is achieved at $l_0= (b/ad)^{1/(d+1)}$. Taking $l=\lfloor l_0 \rfloor$ when $l_0 \ge 1$ and using that $l_0/2 \le l \le l_0$ yields
\bea \label{eq:min.in.l}
\min_{l \in \mathbb{N}_1} \left(al^d+\frac{b}{l} \right)\le a\left(\frac{b}{ad}\right)^{\frac{d}{d+1}} + 2b\left(\frac{ad}{b}\right)^{\frac{1}{d+1}}
= a^{\frac{1}{d+1}}b^{\frac{d}{d+1}} \left(\frac{1}{d^{\frac{d}{d+1}}}+2 d^{\frac{1}{d+1}}\right).
\ena
\bigskip

\noindent {\textbf{Proof of Theorem \ref{multiising}:}}  
By second order stationarity, it suffices to prove the case ${\bf k}={\bf 1}$.
Let $n \ge 2, B_{{\bf 1}}^n$ the block of size $n^d$
as given in \eqref{blockdef}, and $W_{{\bf 1}}^n$ the standardized sum over that block, as in \eqref{wn2}. For any $1 \le l \le n$ write $n=(m-1)l+r, 1 \le r \le l$, and decompose $W_{{\bf 1}}^n$ as the sum of $\xi_{\bf i}/\sqrt{n^dA_n}$ over ${\bf i} \in [m]$, as in \eqref{xidef}.

We apply Theorem \ref{stein:fkg}, handling the two terms on the right hand side of \eqref{stein:fkgbound}. For the first term, using $|X_{{\bf j}}| \le K$, the definition \eqref{xidef} of $\xi_{{\bf i}}$, and the fact that the side lengths of all blocks $D_{\bf i}^l$ are at most $l$, we have
\beas %\label{boundis1}
\left|\frac{\xi_{{\bf i}}}{\sqrt{n^dA_n}} \right| \le B \text{ \ \ with \ } B=\frac{2Kl^d}{\sqrt{n^d A_n}} \text{ \ for all \ } {\bf i} \in [m]^d.
\enas
%\ncolor{Using the bound above, we have 
%\bea \label{1stis}
%\E |\xi_{{\bf i}}|^3 \le B \E \xi_{{\bf i}}^2 
%\ena
%Since $\{ \xi_{{\bf i}}: {\bf i} \in [m]^d \}$ is positively associated and has mean zero, we find that
%\beas
%\sum_{{\bf i} \in [m]^d} \E \xi_{{\bf i}}^2 = \sum_{{\bf i} \in [m]^d} \Var \left(\xi_{{\bf i}}\right) = \Var\left(\sum_{{\bf i} \in [m]^d} \xi_{{\bf i}} \right) - \sum_{{\bf i} \ne {\bf j}} \Cov\left(\xi_{{\bf i}},\xi_{{\bf j}}\right) \le \Var\left(W_{{\bf 1}}^n\right)=1.
%\enas
%Taking the sum in \eqref{1stis} and using the inequality above, we bound the second term as\lcomm{putting second term bound in main theorem}
%\bea \label{boundis2} 
%\sum_{{\bf i} \in [m]^d} \E |\xi_{{\bf i}}|^3 \le B.
%\ena}

Applying Lemma \ref{lemmaising} for the second term, invoking Theorem \ref{stein:fkg} now yields
\beas %\label{eq:d=2bound}
d_1\big({\cal L}(W_{{\bf 1}}),{\cal L}(Z)\big) \le \frac{10Kl^d}{\sqrt{n^d A_n}}+\frac{2\sqrt{2}\kappa_0 \gamma_{\lambda,d}}{\sqrt{\pi}lA_n}.
\enas
Recalling that $A_n>0$ for all $n \in \mathbb{N}_1$, applying the bound \eqref{eq:min.in.l} yields the result, noting that  $l_0=C_{\lambda,\kappa_0,d}^{-1/(d+1)}n^{d/(2d+2)}$ satisfies $1 \le l_0 \le n$ for $n \ge \max \left\{ C_{\lambda,\kappa_0,d}^{2/d},C_{\lambda,\kappa_0,d}^{-2/(d+2)} \right\}$.
%$l=\left\lfloor \left(\frac{\sqrt{2}\kappa_0\left((4\mu_{\lambda}+2\nu_{\lambda})^d-\left(2\nu_{\lambda}\right)^d\right)}{5\lcolor{K}d\sqrt{\pi}\sqrt{A_n}}\right)^{1/(d+1)}n^{d/(2d+2)}\right\rfloor$, using the assumption $A_n > 0$ 
\bbox

\bigskip

To prove Theorem \ref{multiising2} we apply Theorem \ref{stein:fkg2} and use the same techniques as for Theorem \ref{multiising}. We remind the reader that for this result we do not explicitly compute the constants, but index them by the parameters on which they depend.

\bigskip

\noindent {\textbf{Proof of Theorem \ref{multiising2}:}}  We proceed as in the one dimensional case. For $n \ge 2$, and $1 \le l \le n$ we write $n=(m-1)l+r$ with 
$m \ge 1$ and $1 \le r \le l$, and decompose $S_{{\bf k}_q}^n-\E S_{{\bf k}_q}^n$ for $q\in [p]$ as the sum over ${\bf i} \in [m]^d$ of the variables $\xi_{{\bf i},{\bf k}_q}$ given in \eqref{xidef}.

Applying Theorem \ref{stein:fkg2}, we handle the three terms on the right hand side of (\ref{steinbound:fkg2}). For the first term, using the definition \eqref{xidef} of $\xi_{{\bf i},{\bf k}_q}$ and $|X_{\bf t}|\le K$, we have
\beas 
\left|\xi_{{\bf i},{\bf k}_q}\right| \le  B \text{ \ where \ } B=2Kl^d \text{ \ for all \ } {\bf i} \in [m]^d, q \in [p]
\enas
and thus noting that $|\Sigma^{-1/2}|_\infty = n^{-d/2}\psi_{n}$ and that $\Sigma_{jj}=n^dA_n$, we may bound the first term as
\bea \label{multi1st}
\left(\frac{1}{6}+2\sqrt{2}\right)p^3Bn^{-3d/2}\psi_{n}^3\sum_{q=1}^p\Sigma_{q,q} \le \frac{C_{p,K}l^dA_n\psi_{n}^3}{n^{d/2}}.
\ena

For the second term, by Lemma \ref{lemmaising} we have
\begin{multline} \label{beg}
\left(\frac{3}{\sqrt{2}}+\frac{1}{2}\right)p^2n^{-d}\psi_{n}^2 \sum_{q=1}^p \sum_{{\bf i},{\bf j} \in [m]^d,{\bf i} \neq {\bf j}} \E \left(\xi_{{\bf i},{\bf k}_q}\xi_{{\bf j},{\bf k}_q}\right) \\
= C_p\psi_{n}^2 \sum_{q=1}^p \sum_{{\bf i},{\bf j} \in [m]^d,{\bf i} \neq {\bf j}} \E \left(\frac{\xi_{{\bf i},{\bf k}_q}\xi_{{\bf j},{\bf k}_q}}{n^d}\right)
\le \frac{C_{\lambda,\kappa_0,p,d}\psi_{n}^2}{l}.
\end{multline}
Invoking Lemma \ref{lemmaising2} and assumption \eqref{eq:kqks.separated.alphan}
we have
\beas %\label{dest}
\Sigma_{q,s} = \Cov \big(S_{{\bf k}_q}^n,S_{{\bf k}_s}^n \big)  \le \kappa_0 \upsilon_{\lambda}^d \alpha n^d \text{ \ for \ } q \ne s \in [p],
\enas
and hence we may bound the last term as
\begin{multline} \label{dest2}
\left( 2\sqrt{2}p^3Bn^{-3d/2}\psi_{n}^3 + \left(\frac{3}{\sqrt{2}}+\frac{1}{2}\right)p^2n^{-d}\psi_{n}^2 \right) \sum_{q,s \in [p],q\ne s} \Sigma_{q,s} \\
\le \left( \frac{C_{p,K}l^d\psi_{n}^3 }{n^{3d/2}} + \frac{C_p\psi_{n}^2}{n^d}\right)\kappa_0  \upsilon_{\lambda}^d\alpha n^d \\ 
\le \frac{C_{\lambda,\kappa_0d,p,K}\alpha l^d\psi_{n}^3 }{n^{d/2}} + C_{\lambda,\kappa_0,d,p}\alpha \psi_{n}^2.
\end{multline}

By Theorem \ref{stein:fkg2} and \eqref{multi1st}, \eqref{beg} and \eqref{dest2}, 
\beas %\label{multi:d=2bound}
d_{\mathcal{H}_{3,\infty,p}}\big(\mathcal{L}(\Sigma^{-1/2}\mathbf{W}),\mathcal{L}(\mathbf{Z})\big) &\le& \left(C_{p,K}A_n + C_{\lambda,\kappa_0, d,p,K}\alpha \right)\frac{\psi_n^3l^d}{n^{d/2}}  + \frac{C_{\lambda,\kappa_0,d,p}\psi_n^2}{l}   + C_{\lambda,\kappa_0,d,p}\alpha\psi_n^2 \nn \\
&\le& C_{\lambda,\kappa_0,d,p,K}\left(\left(A_n+\alpha \right)\frac{\psi_n^3l^d}{n^{d/2}} + \frac{\psi_n^2}{l}   + \alpha\psi_n^2\right).
\enas
Applying \eqref{eq:min.in.l} to the first two terms in parenthesis, we obtain 
\beas
l_0 = \left(\frac{1}
{d\psi_n (A_n+\alpha)}\right)
^{\frac{1}{d+1}}n^{\frac{d}{2d+2}},
\enas
which satisfies $1 \le l_0 \le n$ for the range of $n$ given in \eqref{eq:thm2rangen}, and applying \eqref{eq:min.in.l} yields the result.
\bbox

\bigskip

We now prove the sufficient condition given above for the invertibility of the covariance matrix of ${\bf S}^n$, and the bound on the norm of its inverse. Recall that a matrix $A \in \mathbb{R}^{p \times p}$ is said to be \textit{strictly diagonally dominant} if $|a_{ii}|-\sum_{j \ne i}|a_{ij}|>0$ for all $i \in [p]$.

\bigskip

\noindent {\textbf{Proof of Lemma \ref{lem:Gershgorin}:}} By Lemma \ref{lemmaising2} and the upper bound on $b$ given in \eqref{eq:cond.b.invert}, for all $q \in [p]$  we have 
\beas
\Sigma_{q,q} - \sum_{1\le s\le p, s\ne q} \left|\Sigma_{q,s}\right| \ge n^{d-1}(nA_n -(p-1)  \kappa_0\upsilon_{\lambda}^d b) > 0.
\enas
Hence, $\Sigma$ is a strictly diagonally dominant matrix, and is therefore invertible by the Gershgorin circle theorem, see for instance Theorem 15.10 of \cite{BR14}.  The final claim of the lemma follows from \cite{AN63}, where it is shown that the bound \eqref{eq:infinity.norm.bound.inverse.lem:Gershgorin} holds for the norm $||C||_{\infty} = \max_i \sum_{j =1}^p |c_{ij}|$, which dominates $|C|_{\infty}$.
\bbox

\subsection{The Ising and percolation models} \label{app2}
In this section we apply Theorems \ref{multiising} and \ref{multiising2} over finite blocks to the net magenetization of the Ising model and the number of points in the supercritical bond percolation model that belong to an infinite cluster. 

\begin{corollary} \label{cor:isingpercolation}
The bounds of Theorem \ref{multiising} (Theorem  \ref{multiising2}) apply to the total magnetization $M_{\bf k}$ of \eqref{mkdef} of the Ising model having infinite volume Gibbs measure in $\mathcal{M}_1$ of \eqref{Ismeasures1}, and to the total number of points $U_{\bf k}$ in \eqref{ukdef} that belong to an infinite cluster in the supercritical bond percolation model in dimensions $d \ge 2$ (when choosing $\alpha$ in \eqref{eq:kqks.separated.alphan} to satisfy \eqref{alphacondition}).
\end{corollary}

\proof We show that the fields in the two models in question posses the properties required for the application of Theorems \ref{multiising} and \ref{multiising2}. First, each must comprise almost surely uniformly bounded variables, must be second order stationary, 
positively associated and must have an exponentially decreasing covariance function. Clearly the uniformly bounded condition holds for both, implying the existence of second moments, so second order stationarity holds under strictly stationary. 

By \cite{ELL06}, any Ising model corresponding to an infinite volume Gibbs measure in $\mathcal{M}_1$ is strictly stationary, and by \cite{new80}, positive association holds whenever the interaction $J$ that appears in the Hamiltonian in \eqref{hamil} is nonnegative.  By \cite{LP68} and \cite{ABF87} (see also \cite{DT15}), the covariance between spins decays at an exponential rate when $\beta > 0, h \ne 0$ and $\beta < \beta_c,h=0$, respectively. 

For the supercritical bond percolation model in dimension 2 or higher, strict stationarity follows from the fact that the vertex connection probability $\theta$ is the same for each bond ${\bf e} \in \mathbb{E}^d$, and that the status of each bond is independent of all others.
Again by \cite{new80}, the binary field $\{ \mathbf{1}_{ \{ |\mathcal{C}({\bf x})|=\infty \} } : {\bf x} \in \mathbb{Z}^d \}$ is positively associated.  By \cite{CCGKS89},  the covariance between $ \mathbf{1}_{ \{ |\mathcal{C}({\bf x})|=\infty \} }$ and $ \mathbf{1}_{ \{ |\mathcal{C}({\bf y})|=\infty \} }$ decays exponentially.

Lastly, for the application of Theorem \ref{multiising2}, assumption \eqref{alphacondition} on $\alpha$ and Lemma \ref{lem:Gershgorin} yield the required invertibility of the covariance matrix.

\bbox

\begin{remark} \label{re:AD15.etc}
The earlier version \cite{GW16} of this work depended on a result of the subsequently withdrawn manuscript \cite{AD15} to handle the Ising model with infinite volume Gibbs measures in $\mathcal{M}_2$ of \eqref{Ismeasures2}. In particular, \cite{GW16} shows how the results of Corollary \ref{cor:isingpercolation} apply to this case of the Ising model under the condition that its covariance function is exponentially decreasing.
\end{remark}

\subsection{ The Voter Model} \label{app3}

In this section, we consider the $d$-dimensional voter model $\{ \eta_t: t \ge 0 \}$ with initial state at time zero given by independent Bernoulli variables with parameter $\theta \in (0,1)$ on each site in $\mathbb{Z}^d$, and the occupation time process $\{T_s^t: t \ge 0\}$, $s\ge 0$, defined in (\ref{occudef}), that records the amount of time in $(s,s+t]$ that the origin spends in state one. 
The work \cite{cox83} demonstrates that $T_0^t$ satisfies the CLT for $d\ge 2$ and any $t>0$. Here we refine the results for $d \ge 7$ in Theorem \ref{voterthm} by applying 
Theorem \ref{stein:fkg} to obtain a rate of convergence in the $L^1$ metric. 

Though the authors are unaware of any results in the literature that yield distributional bounds for the voter model, Stein's method has been successful in obtaining bounds in the anti-voter model, see \cite{RR97}. The anti-voter model differs from the voter model in that after a uniformly neighbor of a vertex is chosen the state of that vertex changes to the opposite state of that neighbor. The work of \cite{RR97} used the  exchangeable pair approach to Stein's method, and positive association was not an ingredient.

By \cite{cox83}, for $d\ge 5$ we have
\bea \label{eq:voterAtdef}
\E T_0^t = \theta t \qmq{and $\lim_{t \rightarrow \infty} A_t =\kappa_2$ for some $\kappa_2 > 0$, where} A_t = \frac{\Var\big(T_0^t\big)}{t}  \qm{for $t>0$.}
\ena
Letting
\bea \label{voterAtdef}
A_s^t = \frac{\Var\big(T_s^t\big)}{t},
\ena
we will show in Lemma \ref{Atbehave} that for all $s \ge 0, \lim_{t\rightarrow \infty}A_s^t = \kappa_2$, and for all $t>0$ that 
$\E T_s^t = \theta t$ and $A_s^t>0$. 

Standardizing $T_s^t$ to have mean zero and variance one we obtain
\bea \label{wtdef}
W_s^t=\frac{T_s^t - \theta t}{\sqrt{A_s^t t}}.
\ena

Now we define the last exit time of the origin in $[0,u]$ and $[0,\infty)$ by
\bea \label{lastexit}
L_u= \sup \{ s \le u : Y_s = {\bf 0} \} \qmq{for $u \ge 0$.} 
\ena
and  
\bea \label{ldef}
L=\sup\{s<\infty: Y_s = {\bf 0}\},
\ena
respectively, where $Y_s$ is a rate 1 simple symmetric random walk in $\mathbb{Z}^d$ starting at the origin at time zero. By definition, it is clear that $L_u$ is increasing in $u$ and $L_u \le L$ for all $u \ge 0$. 
The following theorem provides $L^1$ bounds to the normal for the standardized occupation time $W_s^t$, given in \eqref{wtdef}. Lemma \ref{covbound} shows that the second moment of $L_u$, appearing in the bounds below, is finite.

\begin{theorem} \label{voterthm}
For $W_s^t$ as defined in \eqref{wtdef},  for all $d\ge 7$ and $s \ge 0$,
\begin{multline*}
d_1\big({\cal L}(W_s^t),{\cal L}(Z)\big) \le \left(\frac{180\sqrt{2}\theta(1-\theta) \E L_{2(s+t)}^2}{\sqrt{\pi }(A_s^t)^{3/2}}\right)^{1/2}t^{-1/4} \\
\mbox{for all} \quad t \ge \left(\frac{\sqrt{2}\theta(1-\theta) \E L_{2(s+t)}^2}{5\sqrt{\pi A_s^t}}\right)^{2/3},
\end{multline*}
with $Z \sim {\cal N}(0,1)$, $A_s^t$ is as in \eqref{voterAtdef} and $L_u$ as in \eqref{lastexit}.
\end{theorem}

We also extend Theorem \ref{voterthm} to the multidimensional case, obtaining a bound in Theorem \ref{voterthm2} in 
the metric $d_{\mathcal{H}_{3,\infty,p}}$ to the multivariate normal for ${\bf S}^t = (T_{s_1}^t,\ldots,T_{s_p}^t)$. The results could easily be extended to the case where the occupation times are measured over intervals of varying length, as recorded in the vector
$(T_{s_1}^{t_1},\ldots,T_{s_p}^{t_p})$.

\begin{theorem} \label{voterthm2}
Let $d \ge 7$ and $t>0$. For $p\in \mathbb{N}_1$ let $s_1,\ldots,s_p \ge 0$ be 
such that
\bea \label{eq:kqks.separated.alphan2}
\min_{k,l\in[p],k \ne l}|s_k-s_l| \ge (1-\alpha)t,
\ena
for some $0<\alpha<1$. Let ${\bf S}^t = (T_{s_1}^t,\ldots,T_{s_p}^t)$ where $T_s^t$ is defined as in \eqref{occudef}, assume that the covariance matrix $\Sigma$ of ${\bf S}^t$ is invertible and let 
\beas %\label{eq:defpsin2}
\psi_t = t^{1/2}|\Sigma^{-1/2}|_\infty.
\enas
Then, there exists a constant $C_{p,\theta}$ such that with $\mathbf{Z}$ a standard normal random vector in $\mathbb{R}^p$,
\begin{multline} \label{eq:thm2ranget}
d_{\mathcal{H}_{3,\infty,p}}\big(\mathcal{L}(\Sigma^{-1/2}({\bf S}^t -\theta t)),\mathcal{L}(\mathbf{Z})\big) \le 
C_{p,\theta} \left(\left(\sum_{k=1}^pA_{s_k}^t+\alpha + \frac{1}{t}\right)^{1/2} \psi_t^{5/2} t^{-1/4} + \psi_t^2 \left(\alpha + \frac{1}{t}\right) \right) \\
\mbox{for all} \quad t \ge \left(\psi_t\left(\sum_{k=1}^p A_{s_k}^t +\alpha + \frac{1}{t} \right)\right)^{-2/3}.
\end{multline}
\end{theorem}

From \eqref{voterAtdef} we have $\Var(T_s^t) = t A_s^t$ and we show in Lemma \ref{Atbehave} that $A_s^t \rightarrow \kappa_2$ as $t \rightarrow \infty$. For checking the hypothesis that the covariance matrix of ${\bf S}^t$ is invertible and the quantity $\psi_t$ is of order one, we present the following sufficient condition, proved at the end of this section. In addition, we see from  \eqref{eq:thm2ranget} that if $\psi_t=O(1)$ and $\alpha=O(t^{-\frac{1}{4}})$ then the bound will also be of order $O(t^{-\frac{1}{4}})$. 

\begin{lemma}\label{lem:Gershgorin2}
With $t>0$, $\Sigma$ the covariance matrix of the random vector ${\bf S}^t \in \mathbb{R}^p, p \ge 2$ as given in Theorem \ref{voterthm2} is invertible if for some $b \ge 0$
\bea \label{eq:cond.b.invert2}
\left|s_k-s_l\right| \ge t-b \qmq{for all $k \not =l, k,l \in [p]$, and} b< \frac{t\min_{k\in[p]}A_{s_k}^t-(p-1)\theta(1-\theta)\E[L^2]}{2(p-1)\theta(1-\theta) \E [L]},
\ena
where $A_s^t$ is given in \eqref{voterAtdef} and $L$ in \eqref{ldef}. If $b=\alpha t$ for $\alpha$ satisfying
\beas
0 < \alpha < \min\left\{1,\frac{\min_{k\in[p]}A_{s_k}^t-(p-1)\theta(1-\theta)\E[L^2]/t}{2(p-1)\theta(1-\theta) \E [L]}\right\},
\enas
then the matrix $\Sigma$ is invertible and
\beas
|\Sigma^{-1}|_{\infty} \le \frac{1}{t\left(\min_{k\in[p]}A_{s_k}^t - (p-1)\theta(1-\theta)\left(\E[L^2]/t + 2\alpha \E[L] \right)\right)}.
\enas
\end{lemma}

To prove Theorems \ref{voterthm} and \ref{voterthm2}, we apply the following result implied by (0.8) and (0.9) of \cite{cox83}, shown there using a duality that connects the voter model with a system of coalescing random walks constructed by a time reversal, and tracing the $\{0,1\}$ `opinion' of every site back to its genesis at time zero. 

\begin{lemma}[\cite{cox83}]\label{lem:cox83}
For $t \ge 0$ and $0 \le u \le v$, with $L_u$ as in \eqref{lastexit},
\bea \label{eq:key.cov.identity.voter}
{\rm Cov}(\eta_u({\bf 0}),\eta_v({\bf 0}))=\theta(1-\theta)P\left(L_{u+v} > v-u\right).
\ena
\end{lemma}

We now use Lemma \ref{lem:cox83} to prove the following result which will be used in several places in this section.

\begin{lemma} \label{covbound}
For all $0 \le r < s <t$,
\beas
\int_{r}^{s}\int_{s}^t  \Cov (\eta_u(\mathbf{0}),\eta_v(\mathbf{0})) dv du \le\frac{\theta(1-\theta)}{2} \E[L_{s+t}^2],
\enas
where $L_u$ is as in \eqref{lastexit}. Moreover, with $L$ as in \eqref{ldef} we have $\E[L_{u}^2] \le \E[L^2] < \infty$ for all $u>0$ and $d \ge 7$.
\end{lemma}

\proof  Applying covariance identity \eqref{eq:key.cov.identity.voter} from Lemma \ref{lem:cox83} and using the fact that $L_u$ is increasing in $u$, we obtain 
\beas %\label{covtermvot3}
\int_{r}^{s}\int_{s}^t  \Cov (\eta_u(\mathbf{0}),\eta_v(\mathbf{0})) dv du &=& \theta(1-\theta)\int_{r}^{s}\int_{s}^t P\big(L_{u+v} > v-u\big)dvdu \\
&\le& \theta(1-\theta) \int_{r}^{s}\int_{s}^t P\big(L_{s+t} > v-u\big)dvdu. \nn 
\enas
By a change of variables that preserves the difference $v-u$, we have
\begin{multline*}
\int_{r}^{s}\int_{s}^t P\big(L_{s+t} > v-u\big)dvdu
= \int_0^{s-r}\int_{s-r}^{t-r} P\big(L_{s+t} > v-u\big)dvdu \\ \le \int_0^{s-r}\int_{s-r}^t P\big(L_{s+t} > v-u\big)dvdu.
\end{multline*}
Hence, 
\beas
\int_{r}^{s}\int_{s}^t  \Cov (\eta_u(\mathbf{0}),\eta_v(\mathbf{0})) dv du &\le&  \theta(1-\theta)\int_0^{s-r}\int_{s-r}^t P\big(L_{s+t} > v-u\big)dvdu. \nn \\
&=& \theta(1-\theta)\int_0^{s-r}\int_{s-r-u}^{t-u}P(L_{s+t}>v) dvdu\nn \\
&\le& \theta(1-\theta) \int_0^{s-r}\int_{s-r-u}^{t}P(L_{s+t} >v) dvdu\nn \\
&=& \theta(1-\theta) \int_0^{s-r}\int_u^tP(L_{s+t}>v) dvdu \nn \\
&\le& \theta(1-\theta)\int_0^{t}\int_{u}^{t}P(L_{s+t} >v) dvdu \nn \\
&=& \theta(1-\theta) \int_0^t \int_0^{v}P(L_{s+t} >v) dudv \nonumber \\
&=& \frac{\theta(1-\theta)}{2}\int_0^t 2vP(L_{s+t} >v) dv \nn \\
&\le& \frac{\theta(1-\theta)}{2} \E[L_{s+t}^2] %\label{eq:thetasL2} 
\enas
where change of variables are applied in first and second equality, Fubini's theorem in the third equality, 
and where in the final inequality we apply the standard fact, easily shown using Fubini's theorem, that if $\phi(y)$ is a differentiable function on $[0,\infty)$ satisfying $\phi(0)=0$ then 
%(see, for instance, \cite{Dur10} Lemma 2.2.8), that 
for $Y \ge 0$,
\bea \label{intstdfact}
\E [\phi(Y)] = \int_{0}^{\infty} \phi'(y)P(Y>y)dy.
\ena

Since $L_u \le L$ for all $u\ge 0$, to prove the second claim, it suffices to show that $\E[L^2] < \infty$ for $d \ge 7$. Returning to the integral expression for the second moment, and slightly modifying the calculation for the first moment in \cite{cox83} following (0.11), we have
\bea \label{eq:EL2.1}
\E[L^2] = \int_0^{\infty}2tP(L>t)dt = \int_0^{\infty}2t \int_t^\infty P(Y_s=0)ds \gamma_d dt,
\ena
where $\gamma_d=P(Y_t \not = 0 \,\,\forall t \ge 1),$
known to be positive and increasing for $d \ge 3$. Changing the order of integration in \eqref{eq:EL2.1} yields
\beas 
\E[L^2]
=\gamma_d\int_0^{\infty} \int_0^{s}2t  P(Y_s=\mathbf{0})dt ds 
= \gamma_d\int_0^{\infty}  s^2 P(Y_s=\mathbf{0})ds,
\enas
implying that $\E[L^2] < \infty$ for $d\ge 7$ as, see \cite{cox83} for instance,  $P(Y_s=\mathbf{0})=O(s^{-d/2})$. 

\bbox

The next result shows some important properties of the mean and the variance of $T_s^t$. 

\begin{lemma} \label{Atbehave}
For $d \ge 7$, $t > 0$ and $s \ge 0$, with $T_s^t$ as in \eqref{occudef} and $A_s^t$ as in \eqref{voterAtdef},
\beas
\E T_s^t = \theta t, \text{ \ } A_s^t>0  \text{ \ and \ } \lim_{t \rightarrow \infty} A_s^t =\kappa_2,
\enas
where $\kappa_2$ is given in \eqref{eq:voterAtdef}.
\end{lemma}

\proof
We verify the first claim using the definition of $T_s^t$ and \eqref{eq:voterAtdef} to yield
\beas
\E T_s^t = \E \left[T_0^{t+s} - T_0^s\right] = (t+s)\theta - s\theta = t \theta.
\enas 
The second claim follows from (2.1) in \cite{cox83} and the fact that $P(Y_u = {\bf 0})$ is strictly positive for all $u>0$.    

For the final claim, note that for $s=0$, $A_0^t = A_t$ and the result reduces to \eqref{eq:voterAtdef}. For $s>0$, we have
\begin{multline} \label{varTst}
A_s^t=\frac{\Var(T_s^t)}{t} = \frac{\Var(T_0^{t+s} - T_0^s)}{t} = \frac{\Var(T_0^{t+s}) + \Var(T_0^s) -2 \Cov(T_0^s,T_0^{t+s})}{t}\\
 = \frac{(t+s)A_{t+s}}{t} + \frac{sA_s}{t} - \frac{2 \Cov(T_0^{t+s},T_0^s)}{t}.
\end{multline}  
The first term converges to $\kappa_2$ and the second term tends to zero as $t \rightarrow \infty$. Hence it suffices to show that the last term also tends to zero. 

Writing the covariance in integral form in two parts and applying Lemmas \ref{lem:cox83} and \ref{covbound} to the first and second terms, respectively, we obtain 
\beas
\Cov(T_0^s,T_0^{t+s}) &=& \int_0^s \int_0^s \Cov (\eta_u(\mathbf{0}),\eta_v(\mathbf{0})) dv du + \int_0^s \int_s^{s+t} \Cov (\eta_u(\mathbf{0}),\eta_v(\mathbf{0})) dv du \\
&\le& 2\theta(1-\theta)\int_0^s \int_u^s P(L_{u+v} > v-u) dv du + \frac{\theta(1-\theta)}{2} \E[L_{2s+t}^2] \\
&\le& s^2+\E[L^2],
\enas
and thus the last term in \eqref{varTst} tends to zero as $t \rightarrow \infty$ since $\E[L^2]$ is finite for $d \ge 7$ by Lemma \ref{covbound}.

\bbox

Now we use Lemma \ref{covbound} to prove Lemma \ref{lemmavoter}, which bounds the sum of the covariances between $X_{s,i}^t$ and $X_{s,j}^t$ defined in \eqref{occudef2} over $i\ne j \in [m]$, followed by Lemma \ref{lemmavoter2}, which bounds the covariance between occupation times starting at $r,s$ satisfying $|r-s| \ge t-b$.

\begin{lemma} \label{lemmavoter}
For $t>0$, $s \ge 0$ and $m \in \mathbb{N}_1$, with $X_{s,i}^t$ as in \eqref{occudef2} and $L_u$ as in \eqref{lastexit},
\beas
\sum_{i,j \in [m], i \ne j} \Cov(X_{s,i}^t,X_{s,j}^t) \le \theta(1-\theta)(m-1)\E[L_{2(s+t)}^2].
\enas
\end{lemma}

\proof
Using the definition of $X_{s,i}^t$ in \eqref{occudef2}, we have 
\beas %\label{covtermvot4}
\sum_{i,j \in [m], i \ne j} \Cov(X_{s,i}^t,X_{s,j}^t) &=& 2\sum_{i=1}^{m-1}\sum_{j=i+1}^m \Cov(X_{s,i}^t,X_{s,j}^t)  \\
&=& 2 \sum_{i=1}^{m-1} \Cov \big( X_{s,i}^t,\sum_{j=i+1}^m X_{s,j}^t \big)  \\
&=& 2 \sum_{i=1}^{m-1}  \int_{s+(i-1)t/m}^{s+it/m}\int_{s+it/m}^{s+t}  \Cov (\eta_u(\mathbf{0}),\eta_v(\mathbf{0})) dv du. \\
\enas
Applying Lemma \ref{covbound} to the integrals above, and using the monotonicity property of $L_u$ in $u>0$, we have
\beas
\sum_{i,j \in [m], i \ne j} \Cov(X_{s,i}^t,X_{s,j}^t) \le   \sum_{i=1}^{m-1} \theta (1-\theta)\E[L_{2(s+t)}^2] = \theta(1-\theta)(m-1)\E[L_{2(s+t)}^2].
\enas
\bbox

\begin{lemma} \label{lemmavoter2}
For $t>0$, let $0 \le b \le t$ and let $r,s \ge 0$ satisfy $|r-s| \ge t-b$. Then, with $T_s^t$ as in \eqref{occudef} and $L$ as in \eqref{ldef},
\bea \label{lemmavoter2:ineq}
\Cov(T_s^t,T_r^t) \le  \theta(1-\theta)\left(\E[L^2] + 2b\E[L]\right).
\ena  
\end{lemma}

\proof 
It suffices to consider the case $r \ge s$. First assume that $r < s+t$. Using the definition of $T_s^t$ in \eqref{occudef} and breaking the integral that expresses the covariance we wish to bound into three parts, we have
\begin{multline*} 
\Cov(T_s^t,T_r^t) = \int_s^{s+t} \int_r^{r+t} \Cov (\eta_u(\mathbf{0}),\eta_v(\mathbf{0})) dv du\\
                  = \int_s^{r} \int_r^{r+t} \Cov (\eta_u(\mathbf{0}),\eta_v(\mathbf{0})) dv du 
									+ \int_r^{s+t} \int_{s+t}^{r+t} \Cov (\eta_u(\mathbf{0}),\eta_v(\mathbf{0})) dv du\\
									+ \int_r^{s+t} \int_{r}^{s+t} \Cov (\eta_u(\mathbf{0}),\eta_v(\mathbf{0})) dv du 
\end{multline*} 

Applying Lemma \ref{covbound} to the first two integrals, and using the fact that $L_u \le L$, accounts for the first term in \eqref{lemmavoter2:ineq}. For the last integral, using \eqref{eq:key.cov.identity.voter} from Lemma \ref{lem:cox83}, we obtain
\beas
\int_r^{s+t} \int_{r}^{s+t} \Cov (\eta_u(\mathbf{0}),\eta_v(\mathbf{0})) dv du &=& 2\theta(1-\theta)\int_r^{s+t} \int_{u}^{s+t}P(L_{u+v} > v-u) dv du \\
&\le& 2\theta(1-\theta) \int_r^{s+t} \int_{u}^{s+t}P(L > v-u) dv du \\
&=& 2\theta(1-\theta) \int_r^{s+t} \int_{0}^{s+t-u}P(L > v) dv du \\
&\le& 2\theta(1-\theta) \int_r^{s+t} \int_{0}^{\infty}P(L > v) dv du \\
&=& 2\theta(1-\theta) \int_r^{s+t} \E[L] du \\
&=& 2\theta(1-\theta) (t-r+s) \E[L]  \\
&\le& 2\theta(1-\theta) b \E[L],
\enas
thus accounting for the second term in \eqref{lemmavoter2:ineq}, 
where a change of variables is applied in the second equality, \eqref{intstdfact} in the third equality, and the assumption that $r-s \ge t-b$ in the final inequality.

When $r \ge s+t$, by the positivity of the covariance due to assocation, we have 
\beas
\Cov(T_s^t,T_r^t) &=& \int_s^{s+t} \int_r^{r+t} \Cov (\eta_u(\mathbf{0}),\eta_v(\mathbf{0})) dv du \\
                  &\le& \int_s^{r} \int_r^{r+t} \Cov (\eta_u(\mathbf{0}),\eta_v(\mathbf{0})) dv du.
\enas 
Now applying Lemma \ref{covbound} to the integral above, it is bounded by $\theta(1-\theta) \E [L^2] /2$, hence the claim of the lemma is true in this case as it holds with $b$ replaced by zero.

\bbox

\bigskip

Now we use Theorem \ref{stein:fkg} and Lemmas \ref{covbound} and \ref{lemmavoter} to prove Theorem \ref{voterthm}.

\bigskip

\noindent {\textbf{Proof of Theorem \ref{voterthm}:}} 
Let $m$ be a positive integer and 
\bea \label{xidef3}
\xi_{s,i}^t = \frac{X_{s,i}^t-(t/m)\theta}{\sqrt{A_s^t t}} \qmq{for $i \in [m]$, $s \ge 0$}
\ena
where $X_{s,i}^t$ is defined as in \eqref{occudef2}. We apply Theorem \ref{stein:fkg} to $W_s^t = \sum_{i=1}^m \xi_{s,i}^t$, having mean zero and variance one.
By \cite{EPW67}, the components of the vector $\bsxi_t^s=(\xi_{s,1}^t,\ldots,\xi_{s,m}^t)$ are positively associated as they are increasing functions of the positively associated variables $\{X_{s,i}^t: i \in \mathbb{Z}\}$. 

We now handle the terms on the right hand side of the bound (\ref{stein:fkgbound}). 
For the first term, using the definition of $X_{s,i}^t$ in \eqref{occudef2} we have
\beas
0 \le X_{s,i}^t= \int_{s+(i-1)t/m}^{s+it/m} \eta_u(\mathbf{0})du \le \int_{s+(i-1)t/m}^{s+it/m} 1 ds = \frac{t}{m},
\enas
and thus, by using the definition of $\xi_{s,i}^t$ in \eqref{xidef3},\
\beas %\label{1stvoter}
|\xi_{s,i}^t| \le B \text{ \ where \ } B=\frac{2\sqrt{t}}{m\sqrt{A_s^t}} \text{ \ for all \ } i \in [m].
\enas

Applying Lemma \ref{lemmavoter} for the second term, invoking Theorem \ref{stein:fkg} now yields
\beas %\label{one.dim.voter.tominoverm}
d_1\big({\cal L}(W_s^t),{\cal L}(Z)\big) \le \frac{10\sqrt{t} }{m\sqrt{A_s^t}}+\frac{2\sqrt{2}\theta(1-\theta)m\E[L_{2(s+t)}^2]}{\sqrt{\pi}A_s^t t},
\enas
where $\E[L_{2(s+t)}^2]< \infty$ by Lemma \ref{covbound}. Applying \eqref{eq:min.in.l} with $d=1$ and with $m$ in place of $l$ now yields the result, noting that the lower bound on $t$ in the theorem implies that 
\beas
m_0=\left(\frac{5\sqrt{\pi A_s^t}t^{3/2}}{ \sqrt{2}\theta (1-\theta) \E L_{2(s+t)}^2}\right)^{1/2}
\enas
is at least one. 
\bbox

\bigskip

To prove Theorem \ref{voterthm2} we apply Theorem \ref{stein:fkg2} and use the same techniques as for Theorem \ref{voterthm}. For this result we index constants by the parameters on which they depend, though we do not explicitly compute them.

\bigskip

\noindent {\textbf{Proof of Theorem \ref{voterthm2}:}} 
For $m \in \mathbb{N}_1$ decompose $T_{s_k}^t - \theta t$ for $k \in [p]$ as the sum over $i \in [m]$ of the variables $\xi_i^k = X_{s_k,i}^t - (t/m)\theta$ where $X_{s_k,i}^t$ is given in \eqref{occudef2}.

Applying Theorem \ref{stein:fkg2}, we handle the three terms on the right hand side of (\ref{steinbound:fkg2}). In the calculation below we use that $\E[L]$ and $\E[L^2]$ are finite for $d \ge 7$ by Lemma \ref{covbound}. For the first term, again using the definition of $X_{s_k,i}^t$, we have 
\beas
|\xi_i^k| \le B \text{\ where \ } B = \frac{2t}{m} \text{\ for all $i \in [m]$ and $k \in [p]$.\ }
\enas
Since $|\Sigma^{-1/2}|_\infty = t^{-1/2}\psi_{t}$, and $\Sigma_{k,k}=t A_{s_k}^t$ we may bound the first term as
\bea \label{multvot}
\left(\frac{1}{6}+2\sqrt{2}\right)p^3Bt^{-3/2}\psi_{t}^3\sum_{k=1}^p\Sigma_{k,k} = \frac{C_{p}\psi_t^3 \sqrt{t}}{m}\sum_{k=1}^pA_{s_k}^t.
\ena

For the second term, by Lemma \ref{lemmavoter} we have 
\begin{multline} \label{multvot1}
\left(\frac{3}{\sqrt{2}}+\frac{1}{2}\right)p^2t^{-1}\psi_t^2 \sum_{k=1}^p \sum_{i,j \in [m],i \neq j} \E \left(\xi_i^k \xi_j^k \right) \\
= C_p\psi_t^2 \sum_{k=1}^p \sum_{i,j \in [m],i \neq j} \frac{\Cov \left(X_{s_k,i}^t,X_{s_k,j}^t\right)}{t}
\\ \le
\frac{C_{p,\theta}\psi_t^2m}{t} \sum_{k=1}^p\E[L_{s_k+t}^2]  \le \frac{C_{p,\theta}\psi_t^2 m}{t},
\end{multline}
where we have applied Lemma \ref{covbound} for the final inequality.

Invoking Lemma \ref{lemmavoter2} and assumption \eqref{eq:kqks.separated.alphan2}
we have
\beas %\label{multvot2}
\Sigma_{k,l} = \Cov (T_{s_k}^t,T_{s_l}^t)  \le  \theta(1-\theta)\left(\E[L^2] + 2\alpha t\E[L]\right) \text{ \ for \ } k \ne l \in [p],
\enas
and hence we may bound the last term as
\begin{multline} \label{multvot3}
\left( 2\sqrt{2}p^3Bt^{-3/2}\psi_t^3 + \left(\frac{3}{\sqrt{2}}+\frac{1}{2}\right)p^2t^{-1}\psi_t^2 \right) \sum_{k,l \in [p],k\ne l} \Sigma_{k,l} \\
\le \left( \frac{C_{p}t\psi_t^3 }{mt^{3/2}} + \frac{C_p\psi_t^2}{t}\right)\theta(1-\theta)\left(\E[L^2] + 2\alpha t\E[L]\right) \\ 
\le \frac{C_{p,\theta}  \psi_t^3 }{m\sqrt{t}} + \frac{C_{p,\theta} \psi_t^2}{t } + \frac{C_{p,\theta}\alpha \sqrt{t} \psi_t^3 }{m} + C_{p,\theta}\alpha \psi_t^2.
\end{multline}

By Theorem \ref{stein:fkg2} and \eqref{multvot}, \eqref{multvot1} and \eqref{multvot3}, 
\begin{multline} %\label{multi:d=2bound2}
d_{\mathcal{H}_{3,\infty,p}}\big(\mathcal{L}(\Sigma^{-1/2}(\mathbf{S}^t-\theta t)),\mathcal{L}(\mathbf{Z})\big) 
\le \left(C_p\sqrt{t}\sum_{k=1}^p A_{s_k}^t  + \frac{C_{p,\theta}}{\sqrt{t}} +C_{p,\theta} \alpha \sqrt{t}\right) \frac{\psi_t^3}{m} \\
+ \frac{C_{p,\theta}\psi_t^2 m}{t} + \frac{C_{p,\theta} \psi_t^2}{t} + C_{p,\theta} \alpha \psi_t^2  \\
\le C_{p,\theta}\left(\left(\sqrt{t}\sum_{k=1}^p A_{s_k}^t  + \frac{1}{\sqrt{t}} + \alpha \sqrt{t}\right)\frac{\psi_t^3}{m} + \frac{\psi_t^2m}{t} + \psi_t^2\left(\alpha+\frac{1}{t}\right)\right). \nn
\end{multline}
Applying \eqref{eq:min.in.l} to the first two terms in parenthesis now yields the result, noting that 
\beas
m_0 = \left(\psi_t \left(\sum_{k=1}^p A_{s_k}^t + \alpha + \frac{1}{t}\right)\right)^{\frac{1}{2}}t^{\frac{3}{4}}
\enas
is at least one for the range of $t$ given in \eqref{eq:thm2ranget}.
\bbox

Finally, we present the sufficient condition given above for the invertibility of the covariance matrix of ${\bf S}^t$.

\bigskip

\noindent {\textbf{Proof of Lemma \ref{lem:Gershgorin2}:}}  By Lemma \ref{lemmavoter2} and the upper bound on $b$ given in \eqref{eq:cond.b.invert2}, for all $k \in [p]$  we have
\beas
\Sigma_{k,k} - \sum_{l \in [p], l\ne k} \left|\Sigma_{k,l}\right| \ge tA_{s_k}^t - (p-1)\theta(1-\theta) (\E[L^2]+2b\E[L]) > 0.
\enas
Hence, $\Sigma$ is a strictly diagonally dominant matrix, and the claims follow as in the proof of Lemma \ref{lem:Gershgorin}.
\bbox

\subsection{The contact process} \label{app4}
As detailed in the Introduction, 
in this section we consider the functional 
\bea \label{contactdef2}
D_{s,f}^t = \int_s^{s+t} f(\zeta^{\nu_\lambda}(u)) du \text{ \ for \ } s\ge 0,t > 0,
\ena
where $f$ is a non-constant increasing cylindrical function for the supercritical one dimensional contact process $\{\zeta^{\nu_\lambda}(t):t\ge 0\}$
with infection rate $\lambda$, recovery rate 1, and initial configuration having distribution $\nu_{\lambda}$, the unique non-trivial invariant measure of the process. 
It was proved in \cite{Sch86} that in the supercritical case $D_{0,f}^t$ satisfies the CLT for any cylindrical $f$.
In Theorems \ref{contactthm} and \ref{contactthm2} we provide finite sample bounds in the $L^1$ metric, and a multivariate smooth function metric, for this asymptotic when $f$ is increasing and non-constant by applying Theorems \ref{stein:fkg} and \ref{stein:fkg2}, respectively.

In addition to the time dependent quantity \eqref{contactdef2} studied in this section, 
the results in Section \ref{app1} also allow us to obtain a version of Corollary \ref{cor:isingpercolation} for the number of infected sites at any fixed time $t>0$ in a bounded region for the supercritical multidimensional contact process $\{\zeta_t^{\mathbb{Z}^d}({\bf x}) : {\bf x} \in \mathbb{Z}^d\}$ with intial state having mass one on $\mathbb{Z}^d$. In particular, the values of the process at $t$ are positively associated by Theorem B17 of \cite{Lig99}, and the covariance between $\zeta_t^{\mathbb{Z}^d}({\bf x})$ and $\zeta_t^{\mathbb{Z}^d}({\bf y})$ decays exponentially by Theorem 1.7 of \cite{FV03}. Strictly stationary follows from the fact that at time zero all sites are infected.

For the main goal in this section, we apply Lemma \ref{Sch:lemma2} which provides exponential decay of the covariances of $f(\zeta^{\nu_\lambda}(\cdot))$ in time. We note the brief but important comment before the statement of Lemma 1 of \cite{Sch86}, that the process considered there, and defined in an indirect manner, has the distribution of the contact process started in its unique invariant non-trivial distribution. We recall the definition of $\lambda_*$ given in \eqref{contactcritical}.

\begin{lemma}[Lemma 2 of \cite{Sch86}] \label{Sch:lemma2} 
For $\lambda > \lambda_*$ 
and for any cylindrical $f: \mathscr{P}(\mathbb{Z}) \rightarrow \mathbb{R}$ there exist positive constants  $\gamma=\gamma(f,\lambda)$ and $\kappa=\kappa(f,\lambda)$ depending on $\lambda$ and $f$ such that
\bea \label{cov:contact}
|\Cov (f(\zeta^{\nu_\lambda}(r)),f(\zeta^{\nu_\lambda}(s)))| \le \kappa e^{-\gamma|s-r|} \text{ \ for \ } \{r,s\} \subset (0,\infty).
\ena
\end{lemma}

The results of \cite{Sch86} make use of the fact that any cylindrical function $f$ is a finite linear combination of indicators of the form $I_{\bar{B}}(\cdot)$ given in \eqref{indicatordef} with $|B|< \infty$. Hence, $M_f := |f|_{\infty}<\infty$, implying that the moments of $f$ on the process always exist.

%\ncolor{*** \cite{Sch86} used that 
%\beas
%\E \left[T^{-1} \int_0^T f(\eta^A(t)) dt\right] = \E \left[T^{-1} \int_0^T f(\zeta(t)) dt\right] = \int f d \nu
%\enas
%for some extremal invariant measure $\nu$. However, its reference \cite{Gri79}, only showed
%\beas
%\lim_{T\rightarrow \infty} T^{-1} \int_0^T f(\eta^A(t)) dt = \lim_{T\rightarrow \infty} T^{-1} \int_0^T f(\zeta(t)) dt = \int f d \nu \text{ a.s.}.
%\enas 
%By DCT,\lcomm{what is the dominating function you are using here}\ncomm{The constant $|f|_{\infty} < \infty $. $\E |f|_{\infty} < \infty$} we have 
%\beas
%\lim_{T\rightarrow \infty} \E T^{-1} \int_0^T f(\eta^A(t)) dt = \lim_{T\rightarrow \infty} \E T^{-1} \int_0^T f(\zeta(t)) dt =  \int f d \nu.
%\enas
%***}

The next lemma verifies the positive association required for the application of our main results; let  $Y_{s,i}^{t,m}$ be as in \eqref{yidef}.

\begin{lemma} \label{PAcontact}
For $\lambda > \lambda_*$ and any increasing cylindrical function $f: \mathscr{P}(\mathbb{Z}) \rightarrow \mathbb{R}$, the families $\{f(\zeta^{\nu_\lambda}(t)):t\ge 0\}$, and $\left\{Y_{s,i}^{t,m} : i \in [m]\right\}$ for 
$s \ge 0$, $t > 0$ and $m \in \mathbb{N}_1$, are positively associated.
\end{lemma}
\proof By applying (2.21) in \cite{Lig85}, a corollary of Harris' theorem, the proof of Lemma 1 of \cite{Sch86} demonstrates that $\{\zeta^{\nu_\lambda}(t):t\ge 0\}$ is a positively associated family of random variables. Now both claims follow by the fact that both families in the statement are increasing functions in  $\{\zeta^{\nu_\lambda}(t):t\ge 0\}$.
\bbox  

%\beas 
%Y_{s,i}^{t,m} = \int_{s+(i-1)t/m}^{s+it/m} f(\zeta^{\nu_\lambda}(u))du \text{ \ for \ } i \in \mathbb{Z}.
%\enas
%For simplicity, we drop the subscript $s$ when $s=0$ and the superscript $m$ when $m=t$. Since $f$ is increasing, as shown in the proof of Lemma 1 in \cite{Sch86}, $ \left\{Y_{s,i}^{t,m} : i \in \mathbb{Z}\right\}$ is positively associated. \lcomm{Two items here: a) the proof shows association for the values of the process $f(\zeta^{\nu_\lambda})$, not for integrals of the process, and b) the $Y_k$ variables as defined there on the top of page 1293 are processes in $x \in [0,1]$ for every $k$, note that the time index of the variable $Y_k(x)$ is just $k+x$ for $k \in \mathbb{Z}$, so its not quite lining up with the usage we need. Maybe need to go to the sources in that proof and make a simpler, more direct argument. Also, the association statement should be its own lemma, as its key}

Using the remark above Proposition 1 in \cite{Sch86}, the process  $\{ \zeta^{\nu_\lambda}(t): t \in \mathbb{R} \}$ is strictly stationary, hence we may let
\bea \label{ATFdef}
\frac{\Var(D_{s,f}^t)}{t} = A_{f}^t \text{ \ for all \ } s \ge 0, t>0.
\ena
Since $f$ is non-constant increasing and $\zeta^{\nu_\lambda}(u)$ is strictly stationary, by the Remark under Lemma 1 of \cite{Sch86}, $\Var(f(\zeta^{\nu_\lambda}(u))) > 0$ and does not depend on $u$. Since $f(\zeta^{\nu_\lambda}(u)), u \ge 0$ is positively associated, we have that $A_f^t>0$ for all $t>0$. Using the definition of $D_{s,f}^t$ in \eqref{contactdef2}, we also have $A_f^t \le 2M_f^2 < \infty$.

Standardizing $D_{s,f}^t$ to have mean zero and variance one, we obtain
\bea \label{wTdef}
W_{s,f}^t = \frac{D_{s,f}^t - \E D_{s,f}^t}{\sqrt{tA_{f}^t}}=
\frac{\sum_{i=1}^m \left(Y_{s,i}^{t,m} - \E Y_{s,i}^{t,m} \right)}{\sqrt{tA_{f}^t}}.
\ena
The following Lemma provides a bound on the $L^1$ distance between $W_{s,f}^t$ and the standard normal.

\begin{theorem} \label{contactthm}
Let $\lambda > \lambda_*$, $s \ge 0$ and $t>0$. Then for any non-constant increasing cylindrical function $f: \mathscr{P}(\mathbb{Z}) \rightarrow \mathbb{R}$ , with $W_{s,f}^t$ as in \eqref{wTdef} and $Z \sim {\cal N}(0,1)$,
\beas
d_1\left(\mathcal{L}(W_{s,f}^t),\mathcal{L}(Z)\right) \le \left(\frac{360\sqrt{2}\kappa M_f}{\sqrt{\pi}(A_{f}^t)^{3/2}\gamma^2}\right)^{1/2} t^{-1/4},
\text{ \ for all \ } t \ge \left(\frac{2\sqrt{2}\kappa }{5 \gamma^2 M_f \sqrt{\pi A_{f}^t}}\right)^{2/3}
\enas
where $M_f = |f|_{\infty} <\infty$, and $\kappa$ and $\gamma$ are as in \eqref{cov:contact}.
\end{theorem}

We also prove the following multidimensional version of Theorem \ref{contactthm} for the $p$-vector
\bea \label{eq:contact.vector}
{\bf S}_f^t = (D_{s_1,f}^t,\ldots,D_{s_p,f}^t) \qmq{with $s_1,\ldots,s_p \ge 0$,}
\ena
where $D_{s,f}^t$ is defined in \eqref{contactdef2}.

\begin{theorem} \label{contactthm2}
Let $\lambda>\lambda_*$, $t>0$ and for $p\in \mathbb{N}_1$ let $s_1,\ldots,s_p \ge 0$ satisfy
\beas 
\min_{k,l\in[p],k \ne l}|s_k-s_l| \ge (1-\alpha)t
\enas
for some $0<\alpha<1$, and let $f: \mathscr{P}(\mathbb{Z}) \rightarrow \mathbb{R}$ be a non-constant increasing cylindrical function. Assume that the covariance matrix $\Sigma$ of ${\bf S}_f^t$ in \eqref{eq:contact.vector}
is invertible and let 
\beas %\label{eq:defpsin2}
\psi_f^t = t^{1/2}|\Sigma^{-1/2}|_\infty.
\enas
Then, there exists a constant $C_{p,f,\lambda}$ such that with $\mathbf{Z}$ a standard normal random vector in $\mathbb{R}^p$,
\begin{multline} \label{eq:thm2ranget2}
d_{\mathcal{H}_{3,\infty,p}}\left(\mathcal{L}\left(\Sigma^{-1/2}\left({\bf S}_f^t -\E {\bf S}_f^t \right)\right),\mathcal{L}(\mathbf{Z})\right)   \\
\le C_{p,f,\lambda} \left(\left(A_f^t+\alpha + \frac{1}{t}\right)^{1/2} (\psi_f^t)^{5/2} t^{-1/4} + (\psi_f^t)^2 \left(\alpha + \frac{1}{t}\right) \right) \\
\mbox{for all} \quad t \ge \left(\psi_f^t\left( A_f^t +\alpha + \frac{1}{t} \right)\right)^{-2/3}.
\end{multline}
\end{theorem}

We prove Theorems \ref{contactthm} and \ref{contactthm2} using Theorems \ref{stein:fkg} and \ref{stein:fkg2}. Alternatively, one may apply Theorems \ref{multiising} and \ref{multiising2}{ by breaking $D_{s,f}^t$ into a sum of integrals over intervals of length one and a `remainder integral' over an interval of length $t-\left\lfloor t\right\rfloor$. However, this approach will result in constants larger than the ones stated in Theorems \ref{contactthm} and \ref{contactthm2} which breaks $D_{s,f}^t$ into integrals over intervals of equal, optimal length.

For checking the hypothesis that the covariance matrix of ${\bf S}_f^t$ is invertible and the quantity $\psi_f^t$ is of order one, we present the following sufficient condition. We see from  \eqref{eq:thm2ranget2} that if $\psi_f^t=O(1)$ and $\alpha=O(t^{-\frac{1}{4}})$ then the bound will also be of order $O(t^{-\frac{1}{4}})$.

\begin{lemma} \label{lem:Gershgorin3}
With $t>0$, the covariance matrix $\Sigma$ of the random vector ${\bf S}_f^t \in \mathbb{R}^p, p \ge 2$ in \eqref{eq:contact.vector} is invertible if for some $b\ge 0$ 
\beas
\left|s_k-s_l\right| \ge t-b \qmq{for all $k \not =l, k,l \in [p]$, and} b< 
\left( \frac{tA_f^t\gamma}{2(p-1)\kappa}-\frac{1}{\gamma}\right),
\enas
where $A_f^t$ is given in \eqref{ATFdef} and $\kappa, \gamma$ in \eqref{cov:contact}. If $b=\alpha t$ for some $\alpha$ satisfying
\beas
0 < \alpha < \min\left\{1,\frac{A_f^t\gamma}{2(p-1)\kappa}-\frac{1}{\gamma t}\right\},
\enas
then the matrix $\Sigma$ is invertible and
\beas 
|\Sigma^{-1}|_{\infty} \le \frac{1}{t\left(A_f^t - 2\kappa (p-1)\left(\alpha/\gamma + 1/(\gamma^2t)\right)\right)}.
\enas
\end{lemma}

Since Theorems \ref{contactthm} and \ref{contactthm2} and Lemma \ref{lem:Gershgorin3} have the same conclusions as the voter model results in Theorems \ref{voterthm} and \ref{voterthm2} and Lemma \ref{lem:Gershgorin2}, respectively, only differing by constants depending on the parameters of the model, it suffices to prove `contact process' versions of Lemmas \ref{lemmavoter} and \ref{lemmavoter2} and then the results in this section follow similarly as for their voter model counterparts. Proceeding in this manner, Lemma \ref{lemmacontact} bounds the sum of the covariances between $Y_{s,i}^{t,m}$ and $Y_{s,j}^{t,m}$, defined in \eqref{yidef}, over distinct $i,j \in [m]$, and Lemma \ref{lemmacontact2} bounds the covariance between $D_{r,f}^t$ and $D_{s,f}^t$ where $r$ and $s$ are at least $t-b$ apart for some $0 < b \le t$. 

\begin{lemma} \label{lemmacontact}
For $t>0$, $s \ge 0$ and $m \in \mathbb{N}_1$, with $Y_{s,i}^{t,m}$ as in \eqref{yidef},
\beas
\sum_{i,j \in [m], i \ne j} \Cov\left(Y_{s,i}^{t,m},Y_{s,j}^{t,m}\right) \le \frac{2\kappa m}{\gamma^2}
\enas
where $\kappa$ and $\gamma$ are as in \eqref{cov:contact}. 
\end{lemma}

\proof
Applying Lemma \ref{Sch:lemma2} and using the stationarity of $\zeta^{\nu_\lambda}(\cdot)$, for $1 \le i<j \le m$ we have
\begin{multline*}
\Cov\left(Y_{s,i}^{t,m},Y_{s,j}^{t,m}\right) = \int_{(i-1)t/m}^{it/m}\int_{(j-1)t/m}^{jt/m}  \Cov(f(\zeta^{\nu_\lambda}(u)),f(\zeta^{\nu_\lambda}(v)))du dv \\
\le \kappa\int_{(i-1)t/m}^{it/m}\int_{(j-1)t/m}^{jt/m} e^{-\gamma(u-v)} du dv 
=  \kappa\int_{(j-1)t/m}^{jt/m} e^{-\gamma u } du\int_{(i-1)t/m}^{it/m}e^{\gamma v } dv \\
= \frac{\kappa}{ \gamma^2}\left(e^{\gamma t/2m}-e^{-\gamma t/2m}\right)^2 e^{-(j-i)\gamma t/m}.
\end{multline*}
Summing $\Cov(Y_{s,i}^{t,m},Y_{s,j}^{t,m}) $ over $i \ne j$, we obtain
\beas
\sum_{i,j \in [m],i \ne j} \Cov(Y_{s,i}^{t,m},Y_{s,j}^{t,m}) &=& 2 \sum_{1 \le i < j \le m} \Cov(Y_{s,i}^{t,m},Y_{s,j}^{t,m}) \\
&\le& \frac{2\kappa}{\gamma^2} 
\left(e^{\gamma t/2m}-e^{-\gamma t/2m}\right)^2
\sum_{1 \le i < j \le m}  e^{-(j-i)\gamma t/m}\\
&=&\frac{2\kappa}{\gamma^2} 
\left(e^{\gamma t/2m}-e^{-\gamma t/2m}\right)^2 \sum_{k=1}^{m-1}(m-k)e^{-\gamma k t/m} \\
&=& \frac{2\kappa}{\gamma^2} \left(e^{\gamma t/2m}-e^{-\gamma t/2m} \right)^2 \left(\frac{e^{-\gamma t/m}\left((m-1)-me^{-\gamma t/m}+e^{-\gamma t}\right)}{(1-e^{-\gamma t/m})^2} \right) \\
&=& \frac{2\kappa}{\gamma^2}  \left((m-1)-me^{-\gamma t/m}+e^{-\gamma t}\right) \\
&\le& \frac{2\kappa m}{\gamma^2}, 
\enas
where in the third equality we apply the identity in \eqref{Sum.krk4} with $n=m$ and $w = e^{-\gamma t/m}$.
\bbox

\begin{lemma} \label{lemmacontact2}
For $t>0$, let $0 \le b \le t$ and let $r,s \ge 0$ satisfy $|s-r| \ge t-b$. Then, with $D_{s,f}^t$ as in \eqref{contactdef2},
\beas 
\Cov(D_{r,f}^t,D_{s,f}^t) \le  2\kappa \left(\frac{b}{\gamma} + \frac{1}{\gamma^2}\right),
\enas  
where $\kappa$ and $\gamma$ are as in \eqref{cov:contact}.
\end{lemma}

\proof
It suffices to consider the case $s \ge r$. First assume that $s \le r+t$. Using the definition of $D_{s,f}^t$ in \eqref{contactdef2} and breaking the integral that expresses the covariance we wish to bound into three parts, we have
\begin{multline*} 
\Cov(D_{s,f}^t,D_{r,f}^t) = \int_s^{s+t} \int_r^{r+t} \Cov (f(\zeta^{\nu_\lambda}(u)),f(\zeta^{\nu_\lambda}(v))) du dv\\
                          = \int_s^{s+t} \int_r^{s} \Cov (f(\zeta^{\nu_\lambda}(u)),f(\zeta^{\nu_\lambda}(v))) du dv 
									+ \int_{r+t}^{s+t} \int_{s}^{r+t}\Cov (f(\zeta^{\nu_\lambda}(u)),f(\zeta^{\nu_\lambda}(v))) du dv\\
									+ \int_s^{r+t} \int_{s}^{r+t} \Cov (f(\zeta^{\nu_\lambda}(u)),f(\zeta^{\nu_\lambda}(v))) du dv. 
\end{multline*}
Now applying Lemma \ref{Sch:lemma2} we have
\beas 
\Cov(D_{r,f}^t,D_{s,f}^t) &\le& \kappa \int_s^{s+t} \int_r^{s} e^{-\gamma(v-u)} du dv 
				+ \kappa \int_{r+t}^{s+t} \int_{s}^{r+t} e^{-\gamma(v-u)} du dv \nn \\
			&& \hspace{50pt}+ \kappa \int_s^{r+t} \int_{s}^{v} e^{-\gamma(v-u)} du dv 
									+ \kappa \int_s^{r+t} \int_{v}^{r+t} e^{-\gamma(u-v)} du dv. \nn\\
									&=& \kappa \Bigg(\frac{1-e^{-\gamma(s-r)}-e^{-\gamma t}+e^{-\gamma(s-r+t)}}{\gamma^2} 
				+ \frac{1-e^{-\gamma(r-s+t)}-e^{-\gamma (s-r)}+e^{-\gamma t}}{\gamma^2}  \nn\\
									&& \hspace{80pt} + 2\left(\frac{t-(s-r)}{\gamma} - \frac{1-e^{-\gamma(r-s+t)}}{\gamma^2}\right)\Bigg) \nn\\
							    &=& \kappa \left(\frac{2(t-(s-r))}{\gamma} +\frac{e^{-\gamma(s-r+t)}+e^{-\gamma(r-s+t)}-2e^{-\gamma(s-r)}}{\gamma^2} \right)\nn \\
					&\le& \kappa \left(\frac{2(t-(s-r))}{\gamma} +\frac{2}{\gamma^2} \right)\nn \\
									&\le& 2\kappa \left(\frac{b}{\gamma} + \frac{1}{\gamma^2}\right).
\enas

When $s> r+t$ we have 
\begin{multline*}
\Cov(D_{s,f}^t,D_{r,f}^t) \le \kappa \int_s^{s+t} \int_r^{r+t} e^{-\gamma(v-u)} du dv 
                          = \kappa \int_s^{s+t} e^{-\gamma v} dv  \int_r^{r+t} e^{\gamma u} du \\
                         	= \frac{\kappa}{\gamma^2}\left(e^{-\gamma(s-r-t)}+e^{-\gamma(s-r+t)} -2e^{-\gamma(s-r)}\right) 
									\le \frac{2\kappa}{\gamma^2}.
									\ignore{\ncolor{= \frac{\kappa}{\gamma^2}\left(e^{-\gamma(s-r-t)}\left(1+e^{-2\gamma t} -2e^{-\gamma t}\right)\right)} \\
										\ncolor{=\frac{\kappa}{\gamma^2}\left(e^{-\gamma(s-r-t)}(1-e^{-\gamma t})^2\right) }        
										\le \frac{\kappa}{\gamma^2},}
								\end{multline*}
hence the claim of the lemma is true in this case as it holds with $b$ replaced by zero.
\bbox

The proofs of Theorems \ref{contactthm} and \ref{contactthm2} and Lemma \ref{lem:Gershgorin3} below closely follow those of Theorems \ref{voterthm} and \ref{voterthm2} and Lemma \ref{lem:Gershgorin2} respectively, so only brief outlines are provided. 

\bigskip

\noindent {\textbf{Proof of Theorem \ref{contactthm}:}} 
Following the outline of the proof of Theorem \ref{voterthm}, we apply Theorem \ref{stein:fkg} to the sum of
\beas
\xi_{s,i}^t = \frac{Y_{s,i}^{t,m}- \E Y_{s,i}^{t,m}}{\sqrt{tA_f^t}}, \qmq{which is absolutely bounded by}
B=\frac{2\sqrt{t}M_f}{m\sqrt{A_f^t}}.
\enas
 Applying Lemma \ref{lemmacontact} for the second term of \eqref{stein:fkgbound} yields
\beas
d_1\big({\cal L}(W_{s,f}^t),{\cal L}(Z)\big) \le \frac{10\sqrt{t}M_f }{m\sqrt{A_f^t}}+\frac{4\sqrt{2}\kappa m}{\sqrt{\pi}\gamma^2 tA_f^t}.
\enas
Applying \eqref{eq:min.in.l} with $d=1$ and with $m$ in place of $l$ as in the proof of Theorem \ref{voterthm} now yields the result.

\bbox

\bigskip

\noindent {\textbf{Proof of Theorem \ref{contactthm2}:}} 
We apply Theorem \ref{stein:fkg2} as in the proof of Theorem \ref{voterthm2}. With $\xi_i^k$ defined by $\xi_i^k = Y_{s_k,i}^{t,m} - \E Y_{s_k,i}^{t,m}$, 
we have $|\xi_i^k| \le B$ where $B = 2M_f t/m$. Then bounding the right hand side of \eqref{steinbound:fkg2}, using that $|\Sigma^{-1/2}|_\infty=t^{-1/2}\psi_f^t$ for the first term, applying Lemmas \ref{lemmacontact} and \ref{lemmacontact2} for the second and the last terms respectively and invoking Theorem \ref{stein:fkg2}, we obtain
\begin{multline}
d_{\mathcal{H}_{3,\infty,p}}\big(\mathcal{L}(\Sigma^{-1/2}(\mathbf{S}^t-\E \mathbf{S}^t),\mathcal{L}(\mathbf{Z})\big) 
\le \left(C_{p,f}\sqrt{t} A_f^t  + \frac{C_{p,f,\lambda}}{\sqrt{t}} +C_{p,f,\lambda} \alpha \sqrt{t}\right) \frac{(\psi_f^t)^3}{m} \\
+ \frac{C_{p,f,\lambda}(\psi_f^t)^2 m}{t} + \frac{C_{p,f,\lambda} (\psi_f^t)^2}{t} + C_{p,f,\lambda} \alpha (\psi_f^t)^2  \\
\le C_{p,f,\lambda}\left(\left(\sqrt{t} A_f^t  + \frac{1}{\sqrt{t}} + \alpha \sqrt{t}\right)\frac{(\psi_f^t)^3}{m} + \frac{(\psi_f^t)^2m}{t} + (\psi_f^t)^2\left(\alpha+\frac{1}{t}\right)\right). \nn
\end{multline}
Applying \eqref{eq:min.in.l} to the first two terms in parenthesis now yields the result.

\bbox

Finally to verify Lemma \ref{lem:Gershgorin3}, we follow the same idea as in the proof of Lemma \ref{lem:Gershgorin2} with the help of Lemma \ref{lemmacontact2}.

\bigskip

\section{Proofs of main theorems} \label{pf}
In this section we prove our main results, Theorems \ref{stein:fkg} and \ref{stein:fkg2}, stated in the Introduction. For this purpose we first recall that the pair of random variables $(X,Y)$ is said to be \textit{positive quadrant dependent}, or PQD, if for all $x,y \in \mathbb{R}$ we have
\beas %\label{pqd}
H(x,y) \ge 0 \qmq{where} H(x,y)=P(X>x,Y>y) - P(X>x)P(Y>y).
\enas
It was shown in \cite{Hoe40} (see also Lemma 2 of \cite{Leh66}), that if $(X,Y)$ is PQD and if both $X$ and $Y$ have finite second moments then
\beas %\label{leh:cov}
\Cov \big(X,Y \big) = \int_\mathbb{R} \int_\mathbb{R} H(x,y)dx dy.
\enas
In particular, if $(X,Y)$ is PQD then $\Cov\big(X,Y\big) \ge 0$, and if $\Cov \big(X,Y \big)=0$ then $X$ and $Y$ are independent. 

The following two lemmas are invoked in the proofs of the main theorems. Lemma \ref{new80:lem} is Lemma 3 of \cite{new80}, and allows us to bound $\Cov\big(\phi(X),\psi(Y) \big)$ by a constant times $\Cov\big(X,Y \big)$ for a PQD pair $(X,Y)$ and $\phi$ and $\psi$ sufficiently smooth.

\begin{lemma}[\cite{new80}] \label{new80:lem}
Let the random variables $X,Y$ have finite second moments and be PQD. Then for any real valued functions $\phi$ and $\psi$ that are absolutely continuous on all finite subintervals of $\mathbb{R}$,
\bea \label{Cov.Hxy}
\left| \Cov \big(\phi(X),\psi(Y) \big) \right|
%= \int_\mathbb{R} \int_\mathbb{R} \phi'(x) \psi'(y) H(x,y)dx dy
\le |\phi'|_{\infty} |\psi'|_{\infty} \Cov\big(X,Y \big).
\ena
\end{lemma}

The following result of \cite{new80}, contained in the remark explaining (12) there, provides a version of \eqref{Cov.Hxy} for vector valued functions.
\begin{lemma}[\cite{new80}] \label{new80:remark}
If $\bsxi=(\xi_1,\ldots,\xi_n)$ are positively associated then for all real valued differentiable functions $\phi$ and $\psi$ on $\mathbb{R}^n$,
\beas
\left|\Cov \big(\phi(\bsxi),\psi(\bsxi) \big)\right| \le 3 \sqrt{2} \sum_{i,j=1}^n \left|\frac{\partial \phi}{\partial \xi_i}\right|_{\infty}\left|\frac{\partial \psi}{\partial \xi_j}\right|_{\infty} \Cov \big(\xi_i,\xi_j \big).
\enas
\end{lemma}

As for any fixed $u$ the indicator ${\bf 1}_{X>u}$ is increasing in $X$, when $X$ and $Y$ are positively associated then they are PQD. As positive association is 
preserved by coordinate increasing functions, the following lemma is immediate. 

\begin{lemma} \label{FKG->PQD}
The pair $(X,Y)=(\psi(\bsxi), \phi(\bsxi))$ is PQD whenever $\bsxi$ is positively associated and $\psi(\bsxi)$ and $\phi(\bsxi)$ are coordinate wise increasing functions of $\bsxi$.
\end{lemma}

%\proof For any fixed $x,y, \in \mathbb{R}$ the functions %${\bf 1}_{\psi(\bsxi)>x}$ and ${\bf 1}_{\phi(\bsxi)>y}$ are coordinate wise increasing in $\bsxi$. Hence by (\ref{fkg:def}) we have
%\beas
%0 \le \Cov \left({\bf 1}_{\psi(\bsxi)>x},{\bf 1}_{\phi(\bsxi)>y} \right)= \Cov\left({\bf 1}(X>x),{\bf 1}(Y>y) \right) = H(x,y).
%\enas
%\bbox

Though one form of the $L^1$, or Wasserstein distance $d_1\big({\cal L}(X),{\cal L}(Y)\big)$ is given in (\ref{def.d1.integral}), in the proof that follows we apply the alternative characterization, see \cite{Rac84} for example:
\bea \label{d1:sup.L}
d_1\big({\cal L}(X),{\cal L}(Y)\big) = \sup_{h \in \mathbb{L}} |\E h(X)-\E h(Y)| \qmq{where} \mathbb{L}=\{h: |h(y)-h(x)| \le |y-x|\}.
\ena

\bigskip

\noindent {\bf Proof of Theorem \ref{stein:fkg}}  For given $h \in \mathbb{L}$ let $f$ be the unique bounded solution to the Stein equation
\bea \label{stein:eq}
f'(w)-w f(w) = h(w)-Nh \qmq{where} N h=\E h(Z),
\ena
with ${\cal L}(Z)$ the standard normal distribution. Then, (see e.g. \cite{Che11} Lemma 2.4),
\bea \label{steineq:bounds}
|f'|_{\infty} \le \sqrt{\frac{2}{\pi}} \qmq{and} |f''|_{\infty} \le 2.
\ena
As ${\rm Var}(W)=\sum_{i,j} \sigma_{ij}=1$, we obtain
\beas
\E[f'(W)]&=& \E\left(\sum_{i=1}^m \sigma_i^2 f'(W) + \sum_{i \not = j} \sigma_{ij} f'(W) \right)\\
         &=& \E\left(\sum_{i=1}^m \xi_i^2 f'(W)+ \sum_{i \not = j} \sigma_{ij} f'(W) + \sum_{i=1}^m (\sigma_i^2-\xi_i^2) f'(W) \right).
\enas 
Now letting $W^i=W-\xi_i$, write
\beas
\E[Wf(W)] = \E\sum_{i=1}^m \xi_i f(W) = \E\sum_{i=1}^m \xi_i f(W^i  + \xi_i)
         = \E\sum_{i=1}^m \left[ \xi_i f(W^i) + \xi_i^2 \int_0^1 f'(W^i + u \xi_i) du \right].
\enas
Recalling the Stein equation (\ref{stein:eq}) and subtracting, we obtain
\begin{multline} \label{add.subtract.alt}
\E[h(W)-N h]=\E[f'(W)-Wf(W)]\\
= \E \left( \sum_{i=1}^m \xi_i^2 \left(  \int_0^1 \left(f'(W) -f'(W^i+ u\xi_i)\right)du\right)+  \sum_{i=1}^m (\sigma_i^2 - \xi_i^2) f'(W) \right.\\
\left.  + \sum_{i \not = j} \sigma_{ij}f'(W) - \sum_{i=1}^m  \xi_i f(W^i)\right).
\end{multline}

Using the second inequality in (\ref{steineq:bounds}),  we bound the first term in (\ref{add.subtract.alt}) by
\begin{multline} \label{term:1}
\left| \E \sum_{i=1}^m \xi_i^2 \int_0^1 \left(f'(W) -f'(W^i+ u\xi_i) \right) du \right|\\
= \left| \E \sum_{i=1}^m \xi_i^2 \int_0^1 \int_{u\xi_i}^{\xi_i} f''(W^i+t) dt du  \right|
\le 2 \E \sum_{i=1}^m \xi_i^2 \left(  \int_0^1 \int_{u|\xi_i|}^{|\xi_i|} dt du\right)\\
=  \E \sum_{i=1}^m |\xi_i|^3 \le B \E \sum_{i=1}^m \xi_i^2 \le B,
\end{multline}
using the almost sure bound on the variables $\xi_i$, and that their sum has mean zero and variance 1.

To handle the second term in (\ref{add.subtract.alt}), first note that $W$ and $\xi_i$ are coordinate wise increasing functions of $\bsxi$, and hence PQD by Lemma \ref{FKG->PQD}. Now applying Lemma \ref{new80:lem} with
\beas
g(x)=\left\{  
\begin{array}{cc}
x^2 & |x| \le B\\
B^2 & |x|>B
\end{array}
\right.
\enas
and again using the second inequality in (\ref{steineq:bounds}), we have
\bea
\left| \E \sum_{i=1}^m f'(W)(\sigma_i^2 - \xi_i^2) \right|
&=& \left| \sum_{i=1}^m \Cov \big(f'(W),g(\xi_i) \big) \right|\nonumber \\
&\le& 4B \sum_{i=1}^m   \Cov \big(W,\xi_i \big) = 4B \sum_{i=1}^m \sum_{j=1}^m \sigma_{ij} 
=  4B, %\label{term:2}
\ena
using that $\Var \big(W \big)= \sum_{i,j} \sigma_{ij}=1$.

For the third term in (\ref{add.subtract.alt}), using the nonnegativity of the covariances $\sigma_{ij}$ and the first inequality in (\ref{steineq:bounds}) we obtain
\bea %\label{term:3}
\left| \E\sum_{i \not = j} \sigma_{ij}f'(W) \right| = |\E f'(W)| \sum_{i \not = j} \sigma_{ij}  \le \sqrt{\frac{2}{\pi}} \sum_{i \not = j} \sigma_{ij}.
\ena

For the final term in (\ref{add.subtract.alt}), we note that the variables $W^i$ and $\xi_i$ are coordinate wise increasing in $\bsxi$, hence the pair $(W^i,\xi_i)$ is PQD by Lemma \ref{FKG->PQD}. Applying Lemma \ref{new80:lem} and the first inequality in (\ref{steineq:bounds}) now yields
\bea \label{term:4}
\left| \E \sum_{i=1}^m  \xi_i f(W^i) \right| = \left| \sum_{i=1}^m \Cov\big(\xi_i,f(W^i) \big) \right| \le \sqrt{\frac{2}{\pi}}\sum_{i=1}^m \Cov \big(\xi_i,W^i \big) = \sqrt{\frac{2}{\pi}}\sum_{i \not = j}\sigma_{ij}.
\ena

Summing the bounds (\ref{term:1})-(\ref{term:4}) we find that $|\E h(W)-N h|$ is bounded by the right hand side of (\ref{stein:fkgbound}). Taking supremum over $h \in \mathbb{L}$ and using the characterization of the $d_1$ metric given in (\ref{d1:sup.L}) completes the proof. \bbox

%Alternate ending would be to write the expectation of
%\beas
%\sum_{i \not =j} f'(W)
%\enas

\bigskip

To prove Theorem \ref{stein:fkg2} we apply the following result, which is a small variant of Lemma 2.6 of \cite{Che11}, due to \cite{Bar90}. Let $\mathbf{Z}$ be a standard normal random vector in $\mathbb{R}^p$. For $h: \mathbb{R}^p \rightarrow \mathbb{R}$ let $N h = \E h(\mathbf{Z})$ and for $u \ge 0$ define
\beas
(T_uh)(\mathbf{s})=\E{h(\mathbf{s}e^{-u}+\sqrt{1-e^{-2u}}\mathbf{Z})}.
\enas
We write $D^2 h$ for the Hessian matrix of $h$ when it exists. 

\begin{lemma} \label{multilemma}
For $m \ge 3$ and $h \in L^{\infty}_m(\mathbb{R}^p)$ the function 
\beas
g(\mathbf{s}) = -\int_0^{\infty}[T_uh(\mathbf{s})-N h]du
\enas
solves
\beas
{\rm tr} D^2g(\mathbf{s})-\mathbf{s}\cdot \nabla g(\mathbf{s}) = h(\mathbf{s})-N h,
\enas
and for any $0 \le |{\bf k}|_1 \le m$
\beas
|g^{(\mathbf{k})}|_{\infty}\le \frac{1}{|{\bf k}|_1} |h^{(\mathbf{k})}|_{\infty}.
\enas
Furthermore, for any $\boldsymbol\lambda \in \mathbb{R}^p$ and positive definite $p \times p$ matrix $\Sigma$, $f$ defined by the change of variable
\bea \label{solstein}
f(\mathbf{s}) = g(\Sigma^{-1/2}(\mathbf{s}-\boldsymbol\lambda))
\ena
solves
\bea \label{Stein's eq2}
tr\Sigma D^2 f(\mathbf{s})-(\mathbf{s}-\boldsymbol\lambda) \cdot \nabla f(\mathbf{s}) = h(\Sigma^{-1/2}(\mathbf{s}- \boldsymbol\lambda)) - N h,
\ena
and satisfies
\beas
|f^{(\mathbf{k})}|_{\infty} \le \frac{p^{|\mathbf{k}|_1}}{|\mathbf{k}|_1}|\Sigma^{-1/2}|_{\infty}^{|\mathbf{k}|_1}|h^{(\mathbf{k})}|_{\infty}.
\enas
In particular, if $h \in \mathcal{H}_{m,\infty,p}$ then 
\bea \label{diffbound}
|f^{(\mathbf{k})}|_{\infty} \le \frac{p^{|\mathbf{k}|_1}}{|\mathbf{k}|_1}|\Sigma^{-1/2}|_{\infty}^{|\mathbf{k}|_1} \text{ \ for all \ } 0 \le |\mathbf{k}|_1 \le m.
\ena
\end{lemma}

We apply the same technique as in the univariate case, along with Lemmas \ref{new80:remark} and \ref{multilemma}, to prove our main multivariate theorem.

\bigskip

\noindent {\bf Proof of Theorem \ref{stein:fkg2}} 
Given $h \in \mathcal{H}_{3,\infty,p}$, let $f$ be the solution of (\ref{Stein's eq2}) given by (\ref{solstein}) with $\boldsymbol\lambda={\bf 0}$. Writing out the expressions in (\ref{Stein's eq2}) yields
\bea \label{1st}
\E\left[h(\Sigma^{-1/2}\mathbf{S})- N h \right] = \E\left[\sum_{j=1}^p \sum_{l=1}^p \Sigma_{j,l}\frac{\partial ^2}{\partial s_j\partial s_l}f(\mathbf{S})-\sum_{j=1}^pS_j\frac{\partial}{\partial s_j}f(\mathbf{S})\right] \nn \\
= \E\sum_{j=1}^p\Sigma_{j,j}\frac{\partial ^2}{\partial s_j^2}f(\mathbf{S})+ \E\sum_{j,l \in [p],j \ne l} \Sigma_{j,l}\frac{\partial ^2}{\partial s_j\partial s_l}f(\mathbf{S})-\E\sum_{j=1}^pS_j\frac{\partial}{\partial s_j}f(\mathbf{S}).
\ena
We consider the first term of (\ref{1st}) and deal with each term under the sum separately for $j=1, \ldots,p$. Letting $\sigma^2_{i,j}=\Var\big(\xi_{i,j} \big)$ and $\sigma_{i,j;k,l} = \Cov\big(\xi_{i,j},\xi_{k,l}\big)$, we have
\begin{multline} \label{2nd}
\Sigma_{j,j}\frac{\partial ^2}{\partial s_j^2}f(\mathbf{S}) = \sum_{i=1}^m \sigma_{i,j}^2\frac{\partial ^2}{\partial s_j^2}f(\mathbf{S}) + \sum_{i,k \in [m],i \neq k} \sigma_{i,j;k,j}\frac{\partial ^2}{\partial s_j^2}f(\mathbf{S}) \\ 
= \sum_{i=1}^m \xi_{i,j}^2\frac{\partial ^2}{\partial s_j^2}f(\mathbf{S}) + \sum_{i,k \in [m],i \neq k} \sigma_{i,j;k,j}\frac{\partial ^2}{\partial s_j^2}f(\mathbf{S}) +\sum_{i=1}^m (\sigma_{i,j}^2-\xi_{i,j}^2)\frac{\partial ^2}{\partial s_j^2}f(\mathbf{S}). 
\end{multline}
Now, with $S_{j*i} = S_j-\xi_{i,j}$ we may write the summands of the third term on the right hand side of (\ref{1st}) as
\bea \label{3rd}
S_j \frac{\partial}{\partial s_j}f(\mathbf{S}) 
&=& \sum_{i=1}^m \xi_{i,j} \frac{\partial}{\partial s_j}f(\mathbf{S}) \nonumber \\
&=& \sum_{i=1}^m \xi_{i,j}\frac{\partial}{\partial s_j}f(S_1,\ldots, S_{j*i} , \ldots, S_p) \nonumber \\ 
&\text{ \ \ \ }& + \sum_{i=1}^m \xi_{i,j}^2 \int_0^1 \frac{\partial ^2}{\partial s_j^2} f(S_1,\ldots, S_{j*i}+u\xi_{i,j} , \ldots ,S_p) du .
\ena

Substituting (\ref{2nd}) and (\ref{3rd}) into (\ref{1st}), we obtain
\bea \label{4th}
&&\E\left[h(\Sigma^{-1/2}\mathbf{S}) - N h \right] \nn\\
&&\hspace{20pt}= \E\sum_{j=1}^p \sum_{i=1}^m \xi_{i,j}^2 \int_0^1 \left( \frac{\partial ^2}{\partial s_j^2}f(\mathbf{S}) -\frac{\partial ^2}{\partial s_j^2} f(S_1,\ldots, S_{j*i}+u\xi_{i,j} , \ldots ,S_p) \right) du \nonumber \\
&&\hspace{30pt}  + \E\sum_{j=1}^p \sum_{i=1}^m (\sigma_{i,j}^2 - \xi_{i,j}^2) \frac{\partial ^2}{\partial s_j^2}f(\mathbf{S}) - \E \sum_{j=1}^p\sum_{i=1}^m \xi_{i,j}\frac{\partial}{\partial s_j}f(S_1,\ldots, S_{j*i} , \ldots, S_p) \nonumber \\
&&\hspace{30pt} + \E\sum_{j=1}^p\sum_{i,k \in [m],i \neq k} \sigma_{i,j;k,j}\frac{\partial ^2}{\partial s_j^2}f(\mathbf{S}) + \E\sum_{j,l \in [p],j \ne l} \Sigma_{j,l}\frac{\partial ^2}{\partial s_j\partial s_l}f(\mathbf{S}).
\ena

We handle these five terms separately. For the first term in (\ref{4th}), using \eqref{diffbound} we have
\bea \label{5th}
\Bigg|\E\sum_{j=1}^p \sum_{i=1}^m &\xi_{i,j}^2& \int_0^1 \left( \frac{\partial ^2}{\partial s_j^2}f(\mathbf{S}) -\frac{\partial ^2}{\partial s_j^2} f(S_1,\ldots, S_{j*i}+u\xi_{i,j} , \ldots ,S_p) \right) du\Bigg| \nonumber \\
&=& \left|\E\sum_{j=1}^p  \sum_{i=1}^m \xi_{i,j}^2 \int_0^1 \int_{u\xi_{i,j}}^{\xi_{i,j}}\frac{\partial ^3}{\partial s_j^3} f(S_1, \ldots , S_{j*i}+t, \ldots,S_p)dtdu\right| \nonumber \\
&\le& \frac{p^3}{3}|\Sigma^{-1/2}|_{\infty}^3\E\sum_{j=1}^p  \sum_{i=1}^m \xi_{i,j}^2 \int_0^1 \int_{u|\xi_{i,j}|}^{|\xi_{i,j}|}dtdu \nonumber \\
&=& \frac{p^3}{6}|\Sigma^{-1/2}|_{\infty}^3\sum_{j=1}^p  \sum_{i=1}^m \E|\xi_{i,j}|^3  \le \frac{p^3}{6}|\Sigma^{-1/2}|_{\infty}^3B\sum_{j=1}^p  \sum_{i=1}^m \E\xi_{i,j}^2 \nn \\
&\le& \frac{p^3}{6}|\Sigma^{-1/2}|_{\infty}^3 B \sum_{j=1}^p\Sigma_{j,j},
\ena
using the almost sure bound on the variables $\xi_{i,j}$, and that their sum $S_j$ over $i$ from $1$ to $m$ has mean zero.

For the second term in (\ref{4th}), we have
\begin{multline*} %\label{eq:conditioning.attempt}
\left|\E\sum_{j=1}^p \sum_{i=1}^m (\sigma_{i,j}^2 - \xi_{i,j}^2) \frac{\partial ^2}{\partial s_j^2}f(\mathbf{S}) \right| = \left|\sum_{j=1}^p \sum_{i=1}^m  \Cov\left(\frac{\partial ^2}{\partial s_j^2}f(\mathbf{S}),\xi_{i,j}^2\right)\right|\\
=\left|\sum_{j=1}^p \sum_{i=1}^m  \Cov\left(\frac{\partial ^2}{\partial s_j^2}f(\mathbf{S},\xi_{i,j}),g({\bf S},\xi_{i,j})\right)\right|,
\end{multline*}
where with some abuse of notation we let $f({\bf s},x)=f({\bf s})$, and define
\beas
g({\bf s},x)=\left\{  
\begin{array}{cc}
x^2 & |x| \le B\\
B^2 & |x|>B,
\end{array} \right.
\enas
for all ${\bf s} \in \mathbb{R}^p$ and $x \in \mathbb{R}$. Applying Lemma \ref{new80:remark}
and using the bound \eqref{diffbound} and the fact that $({\bf S},\xi_{i,j})$ are positively associated for all $i,j$, we obtain
\bea %\label{6th}
\left|\sum_{j=1}^p \sum_{i=1}^m  \Cov\left(\frac{\partial ^2}{\partial s_j^2}f(\mathbf{S}),\xi_{i,j}^2\right)\right| 
&\le& 3\sqrt{2} \left|\sum_{j=1}^p \sum_{i=1}^m \sum_{l=1}^p  \left|\frac{\partial^3}{\partial s_l\partial s_j^2} f \right|_{\infty}\left|\frac{\partial g}{\partial x}\right|_{\infty} \Cov\left(S_l,\xi_{i,j}\right)\right| \nonumber \\
&\le& 2\sqrt{2}p^3|\Sigma^{-1/2}|_{\infty}^3B\sum_{j=1}^p \sum_{i=1}^m \sum_{l=1}^p  \Cov\big(S_l,\xi_{i,j}\big) \nonumber \\
&=&  2\sqrt{2}p^3|\Sigma^{-1/2}|_{\infty}^3B\sum_{j,l=1}^p \sum_{i,k=1}^m  \sigma_{i,j;k,l} \nonumber \\
&=&  2\sqrt{2}p^3|\Sigma^{-1/2}|_{\infty}^3B \sum_{j,l=1}^p \Sigma_{j,l}.
\ena

For the third term in (\ref{4th}), again applying Lemma \ref{new80:remark} and arguing as for the second term, we have
\bea %\label{7th}
&&\left|\E\sum_{j=1}^p \sum_{i=1}^m  \xi_{i,j}\frac{\partial}{\partial s_j}f(S_1,\ldots, S_{j*i} , \ldots, S_p)\right| \nn\\
&&\hspace{20pt} = \left|\sum_{j=1}^p \sum_{i=1}^m \Cov\left(\xi_{i,j},\frac{\partial}{\partial s_j}f(S_1,\ldots, S_{j*i} , \ldots, S_p) \right)\right| \nonumber \\
&&\hspace{20pt} \le 3\sqrt{2}\left|\sum_{j=1}^p \sum_{i=1}^m \left(\sum_{l \in [p]\setminus \{j\}}  \left|\frac{\partial^2}{\partial s_l\partial s_j}f\right|_{\infty}  \Cov\left(\xi_{i,j},S_l\right) + \left|\frac{\partial^2}{\partial s_j^2}f\right|_{\infty}\Cov\left(\xi_{i,j},S_{j*i} \right) \right)\right| \nonumber \\
&&\hspace{20pt} \le \frac{3}{\sqrt{2}}p^2|\Sigma^{-1/2}|_{\infty}^2\left( \sum_{j,l \in [p],j \ne l}\sum_{i=1}^m\Cov\left(\xi_{i,j},S_l\right)  + \sum_{j=1}^p \sum_{i=1}^m \Cov \left(\xi_{i,j},S_{j*i} \right)\right) \nonumber \\
&&\hspace{20pt} = \frac{3}{\sqrt{2}}p^2|\Sigma^{-1/2}|_{\infty}^2 \left(\sum_{j,l \in [p],j \ne l}\sum_{i,k=1}^m \sigma_{i,j;k,l} +  \sum_{j=1}^p\sum_{i,k \in [m],i \neq k} \sigma_{i,j;k,j} \right) \nonumber \\
&&\hspace{20pt} = \frac{3}{\sqrt{2}}p^2|\Sigma^{-1/2}|_{\infty}^2 \left(\sum_{j,l \in [p],j \ne l}\Sigma_{j,l} +  \sum_{j=1}^p\sum_{i,k \in [m],i \neq k} \sigma_{i,j;k,j} \right).
\ena

For the fourth and the fifth terms in (\ref{4th}), again using \eqref{diffbound} we have
\bea %\label{8th}
\left|\sum_{j=1}^p\sum_{i,k \in [m],i \neq k} \sigma_{i,j;k,j}\frac{\partial ^2}{\partial s_j^2}f(\mathbf{S})\right| \le \frac{p^2}{2}|\Sigma^{-1/2}|_{\infty}^2\sum_{j=1}^p\sum_{i,k \in [m],i \neq k} \sigma_{i,j;k,j}
\ena
and
\bea \label{9th}
\left|\E\sum_{j,l \in [p],j \ne l} \Sigma_{j,l}\frac{\partial ^2}{\partial s_j\partial s_l}f(\mathbf{S})\right| \le  \frac{p^2}{2}|\Sigma^{-1/2}|_{\infty}^2\sum_{j,l \in [p],j \ne l} \Sigma_{j,l}.
\ena
Summing the bounds (\ref{5th})-(\ref{9th}) we find that $\left|\E\left[h(\Sigma^{-1/2}\mathbf{S}) - N h \right]\right|$ is bounded by the right hand side of \eqref{steinbound:fkg2}. Taking supremum over $h \in \mathcal{H}_{3,\infty,p}$ and using the definition \eqref{smoothdef} of $d_{\mathcal{H}_{m,\infty,p}}$  completes the proof.\bbox

\def\cprime{$'$}

\end{document}